\documentclass[11pt, oneside, reqno]{amsart}
\usepackage[shortlabels, inline]{enumitem}
\usepackage{mathtools, mathrsfs, xargs, amssymb}
\usepackage[T1]{fontenc}
\usepackage{libertine}
\usepackage[top=1.2in, bottom=1.2in, left=1.2in, right=1.2in]{geometry}
\usepackage{stmaryrd}

\usepackage[dvipsnames]{xcolor}
\usepackage[
  colorlinks=true,
  linkcolor=Maroon,
  citecolor=Maroon,
  urlcolor=Maroon
]{hyperref}

\newcommand{\fin}{
  \bibliographystyle{amsalpha} 
  \bibliography{refs}
  \end{document}
}

\makeatletter
\renewcommand\subsubsection{\@startsection{subsubsection}{3}{\z@}%
  {.2\linespacing\@plus.2\linespacing}{-.5em}{\normalfont\itshape}}
\makeatother

\usepackage{defs}

\NewIndexRange{N}{N}{N\geq 2}

\newtheorem{theorem}{Theorem}[section]
\newtheorem{proposition}[theorem]{Proposition}
\newtheorem{lemma}[theorem]{Lemma}
\newtheorem{corollary}[theorem]{Corollary}
\theoremstyle{definition}
\newtheorem{definition}[theorem]{Definition}
\newtheorem{assumption}[theorem]{Assumption}
\theoremstyle{remark}
\newtheorem{remark}[theorem]{Remark}
\numberwithin{equation}{section}

\makeatletter
\AddToHook{env/enumerate/before}{\if@inlabel\leavevmode\fi}
\AddToHook{env/itemize/before}{\if@inlabel\leavevmode\fi}
\makeatother

\begin{document}

\title[Interacting Brownian Particles with Local Time Driven Jumps and Memory]{
  Interacting Brownian Particles with Local Time Driven Jumps and Memory}

\author{Flavio Argentieri}
\address{Flavio Argentieri, Department of Mathematics, The
University of Texas at Austin}
\email{flavio.argentieri@utexas.edu}

\author{Gordan {\v Z}itkovi{\' c}}
\address{Gordan {\v Z}itkovi{\' c}, Department of Mathematics, The
University of Texas at Austin}
\email{gordanz@math.utexas.edu}

\begin{abstract} We study a system of Brownian particles on the
  real line whose jump intensities are excited by pairwise
  collision local times and decay exponentially between
  encounters. For every finite number of particles and every
  interaction strength, we construct a global solution and prove
  uniqueness in law. Under mean-field scaling, exchangeable
  chaotic initial data, and uniform integrability of the initial
  intensities, we prove propagation of chaos on path space. The
  limit is a nonlinear jump-diffusion
  whose intensity is driven by the density of its own position; it is
  unique among such processes with locally bounded mean intensity. The proof combines uniform
  intensity and local-time estimates, tightness of particle and
  pair empirical measures, identification of conditional
  collision local times, and a coupling argument for nonlinear
  uniqueness. No moment assumptions are imposed on the initial
  positions or the jump-size distribution. 
\end{abstract}


\thanks{
  During the preparation of this work both authors were
  supported by the National Science Foundation
  under Grant DMS-2307729. Any opinions, findings and
  conclusions or recommendations expressed in this material are
  those of the authors and do not necessarily reflect the views of the
  National Science Foundation (NSF).
}

\maketitle

\section{Introduction}\label{s-intro}
We study a Brownian particle system on the real line in which
particles jump at rates that are excited by pairwise contacts and
relax exponentially in between, so that each rate carries a fading
memory of past encounters.
More precisely,
contacts are measured by collision local times,
and jump sizes follow a fixed probability 
  measure $\nu$ on
$E = \R \setminus \set{0}$,
so that the
$N$-particle system with the interaction-strength parameter
$\gm_N \geq 0$ follows
\begin{align*}
  \diffd X^i_t = \diffd B^i_t + \int_{E} z \, \mu^i(\dt,\dz),\ \
  \diffd\ld^i_t = -\ld^i_t \, \dt + \gm_N \textstyle\sum_{j \neq i}
  \diffd\ell^{i,j}_t,\ 
  i = 1, \dots, N,
\end{align*}
where $B^1, \dots, B^N$ are independent Brownian motions, $\mu^i$ is
the jump measure of particle $i$, with compensator
$\ld^i_t \, \dt \, \nu(\dz)$, and $\ell^{i,j}$ is the collision local
time of particles $i$ and $j$. 
We work in one dimension because in dimension two or higher
independent Brownian particles almost surely never meet, so collision
local times vanish; a higher-dimensional version would need a
mollified interaction, as in the moderate-interaction models discussed
below.

For the finite system
we construct a global solution, prove uniqueness in law, and show that
the solution is a measurable function of the initial data, the
Brownian motions and the Poisson random measures that drive the
jumps.
We then consider the mean-field
scaling $N\gm_N\to\kp\geq0$. Our main result is qualitative
propagation of chaos on path space for exchangeable, chaotic
initial conditions with uniform integrability of the initial
intensities. No moment assumptions are imposed on the initial
positions or on $\nu$, and the initial positions are not required
to have a density.
We also show that the limit is the unique solution, among
solutions with locally bounded mean intensity (Definition~\ref{d-nl}),
of the nonlinear
system
\begin{align*}
  \diffd X_t = \diffd B_t + \int_E z \, \mu(\dt,\dz), \ \ 
  \diffd\ld_t = \bigl(-\ld_t + \kp \, u(t,X_t)\bigr) \, \dt,
\end{align*}
where $\mu$ has compensator $\ld_t \, \dt \, \nu(\dz)$ and $u(t,\cdot)$
is the probability density of the position $X_t$. This density
exists for every $t > 0$, whether or not $X_0$ has a density,  and satisfies
$\norm{u(t,\cdot)}_\infty \leq C_T \, t^{-1/2}$
for $0 < t \leq T$ (Lemma~\ref{l-densbd}), so the forcing term
$\kp \, u(t,X_t)$ is integrable near $t = 0$.

The collision local times thus become a
density-dependent excitation in the limit, and the intensity carries
an exponentially decaying memory of the density the particle has met
along its path. Because only the interaction strength is scaled, the
jumps do not disappear since a tagged particle keeps jumping, and the jump
noise cancels only after averaging over the particles, as it must,
since the empirical measures have a deterministic limit
(Remark~\ref{r-jumps}). We also prove convergence of the pair
empirical measures, which carry the collision local times, as
Sznitman did for the Burgers model \cite{Szn86}.

The proof has to deal with three problems which, as
far as we know, do not occur together in earlier work.
First, the
excitation is singular, since it is carried by local times, and could
in principle lead to an accumulation of jumps; the unbounded
intensities must also be controlled uniformly in $N$, under uniform
integrability of the initial intensities alone. Second, the collision
local time is not a continuous functional of the paths in the
Skorokhod $J_1$ topology, so its convergence does not follow from
tightness of the particle paths. Third, the limiting intensity depends
on the density of the law at the particle's own position, which is not
a continuous functional of the law in the weak topology.
Section~\ref{ss-heuristic} outlines how these difficulties are
overcome.

\medskip

Our work is part of the large literature on mean-field
limits of interacting particle systems, which goes back to Kac's
notion of chaos and its propagation \cite{Kac56}. Probabilistic
approaches to the Boltzmann equation were developed by Tanaka
\cite{Tan78} and Sznitman \cite{Szn84a}; Vlasov's equation for plasmas
\cite{Vla68}, where particles interact through long-range forces, is
another source. General accounts of mean-field theory and propagation of
chaos, with additional references, can be found in
\cite{Szn91,Mel96,Gol16,JabWan17,CarDel18a,CarDel18b,ChaDie22a,ChaDie22b}.

McKean's work on nonlinear Markov
processes \cite{McK66,McK67} carried these ideas beyond kinetic
theory. A standard method for proving propagation
of chaos consists of three steps: establishing tightness of
the empirical measures, identifying their possible limit points as
solutions of a nonlinear problem, and proving uniqueness of the
limiting law. The compactness and nonlinear martingale-problem methods
behind this strategy were further developed in
\cite{Oel84,Fun84,Leo86,Gar88,Szn91} and many other papers.
These compactness arguments are qualitative; quantitative approaches
to propagation of chaos in kinetic and diffusive settings are
developed in \cite{MisMou13} and \cite{Lac23}, respectively.

A model that is particularly relevant to our work was proposed by
McKean \cite{McK67} and studied by Sznitman \cite{Szn86} and by
Sznitman and Varadhan \cite{SznVar86}. In that model, the particles'
drifts are affected by their pairwise collision local times, leading
to the viscous Burgers equation as the limiting Fokker--Planck
equation. Our model uses a related local-time averaging mechanism, but
here the local times act
on the jump intensities, and only through them on the positions. Our identification of the
limit partly follows Sznitman's approach \cite{Szn86}. More precisely, 
for two
independent copies of the limit, the conditional expectation of their
collision local time given the first copy, $X^1$, is
$\int_0^t u(s, X^1_s) \, \ds$, and this identifies the averaged
collision input, although
here paths are discontinuous and each jump restarts a Brownian piece of
the pairwise distances.

Local interactions also arise in other settings, such as
moderate-interaction models  \cite{Oel85,JouMel98}, where suitably rescaled
interaction kernels with shrinking spatial range lead to nonlinear
equations with local density dependence.  More generally, our
limit belongs to the class of McKean--Vlasov equations whose
coefficients depend on the density of the law at the
particle's own position, a dependence that is not
continuous in the weak topology; see, e.g.,
\cite{BarRoc18,BarRoc20}.  Local times also drive the
dynamics of ranked particle systems, as in the Atlas model
\cite{BanFerKar05,Shk12}, and of self-repelling diffusions
\cite{NorRogWil87}. More recently, Baker, Hambly and Jettkant
\cite{BakHamJet25} studied reflected Brownian particles with interactions
driven by boundary local times, proving a mean-field limit in the symmetric
case up to a possible breakdown time. Jettkant \cite{Jet25} treated the
Atlas model through boundary-interacting McKean--Vlasov equations.
Interactions through hitting times appear in probabilistic formulations of
the supercooled Stefan problem \cite{DelNadShk22}.

\medskip

Jumps, too, already appeared in kinetic models of propagation of chaos.
General treatments of mean-field systems with jumps include Graham
\cite{Gra92a,Gra92b}, Graham and M\'el\'eard \cite{GraMel97}, and,
more recently, Jourdain, M\'el\'eard and Woyczynski \cite{JouMelWoy08}
and Andreis, Dai Pra and Fischer \cite{AndDaiFis18}. To the best of
our knowledge, however, 
 mean-field
limits of systems in which jump intensities are driven by collision
local times, and carry a fading memory of past encounters, have not
been studied before. 

We emphasize that our results are qualitative. While quantitative rates of propagation of
chaos and the study of fluctuations are natural next questions, entropy
methods such as those of \cite{Lac23} do not seem to apply directly,
since collision local times are not regular functionals of the
empirical measure.

\subsection{Applications}\label{ss-app} \label{ss-appl-react}
Our main motivation is the movement of proteins along
DNA, a one-dimensional medium in which particle dynamics depend on the
history of pairwise encounters. The mechanism itself, a combination of diffusion,
encounter-dependent excitation with fading memory, and jumps,
also appears in other contexts (sometimes in higher dimensions
where it motivates the structure of the model rather than its
one-dimensional form): \emph{diffusion-mediated reactions}, with
contact-dependent reactivity \cite{Gre20,BenGre22};
\emph{biological dispersal}, with rapid relocation represented by
position jumps \cite{OthDunAlt88}; and \emph{cell motility}, with
movement regulated by a relaxing internal chemical signal
\cite{RidDalPeaCha23}.

One-dimensional  diffusion with occasional
long jumps is a standard description of how proteins search for
their targets on DNA. A DNA-binding
protein alternates between sliding along the DNA, a one-dimensional
random walk that has been observed directly in living cells
\cite{HamLerMah12}, and episodes in which it detaches, diffuses
through the surrounding medium and rebinds elsewhere
\cite{BerWinvHi81,SluMir04}. Along the DNA, these episodes appear as
jumps. Indeed, 
since DNA is folded in space, a short
excursion through the medium can carry the protein far along the
chain, so the jumps can be long \cite{LomAmbMet05}, which is one
reason to avoid moment assumptions on $\nu$. Moreover, DNA-bound
proteins detach faster when other proteins are present
\cite{GraJohMar11}, so contacts themselves can trigger
relocation.

If a contact changes a
protein's internal state, for instance its conformation or how
strongly it binds to the DNA, the effect of the contact lasts until
that state relaxes, much like the internal chemical signal in the
motility model of \cite{RidDalPeaCha23}. In our model each contact
raises the particle's propensity to relocate, and this effect wears
off gradually after the particles separate. Particles at the same
location can therefore have different jump rates, because their recent
encounters differ, and a particle can stay active after it has left a
crowded region. Memory of this kind is absent
from the encounter-based approach to diffusion-mediated reactions
\cite{Gre20,BenGre22,Gre22}, in which the accumulated local time
on a fixed target controls reactivity but never fades. A natural
extension, not pursued here, lets jumps act on discrete internal
states instead of positions, as in \cite{Bre22}, so that priming
affects the outcome of later encounters rather than the location
of the particle.

\subsection{Outline}\label{ss-outline}
Section~\ref{s-model} presents the model, the main theorem, and some
heuristics, and then establishes existence,
non-explosion, and uniqueness in law for the finite system.
Section~\ref{s-bounds} derives uniform superlinear moment and increment 
bounds, and uses them to prove tightness of the particle and pair empirical
measures. Section~\ref{s-nonlin} studies the nonlinear system,
identifies its density and conditional collision local time, and
proves uniqueness of its law. Section~\ref{s-ident} identifies all limit
points and completes the proof of propagation of chaos.
Appendix~\ref{a-construction} contains the proofs of the construction
and of the Poisson representations for the finite system.
Appendix~\ref{a-ui-function} constructs the uniform-integrability
function used in Section~\ref{ss-bset}, and
Appendix~\ref{a-localtimes} collects the local-time estimates used in
the compactness and identification arguments.
Appendix~\ref{a-nonlinear-representation} proves the Poisson
representation of the nonlinear law, and
Appendix~\ref{a-measure-tools} gives the proofs of
Lemmas~\ref{l-ext} and~\ref{l-pass}.

\subsection{Notation and conventions}\label{ss-notation}
For $N \in \N$ we write $[N] := \set*{1, \dots, N}$. We write
$a \wedge b := \min(a,b)$ and set $\inf \emptyset := \infty$. The sign
function is $\sgn(x) := x/\abs*{x}$ for $x \neq 0$, with
$\sgn(0) := 0$. For a metric space $S$, $\sP(S)$ is the set of Borel
probability measures on $S$ with the topology of weak convergence,
and $\delta_x$ is the Dirac mass at $x \in S$. The image of a measure
$m$ under a map $f$ is $\pushf{f}{m}$. If $S$ is also locally
compact, $C_c(S)$ is the space of continuous functions on $S$ with
compact support. The space $D([0,\infty),S)$ of \cd\ paths carries
the Skorokhod $J_1$ topology, and the space $C([0,\infty),S)$ of
continuous paths carries the locally uniform topology;
$C_\uparrow([0,\infty)) \subseteq C([0,\infty),[0,\infty))$ is the
closed subset of nondecreasing paths starting at $0$. For a locally
compact Polish space $S$, $\sM_p(S)$ is the space of locally finite
integer-valued measures on $S$ with the vague topology, generated by
the maps $m \mapsto \int f \, \diffd m$, $f \in C_c(S)$. Integrals without an indicated
domain are over the whole space.

We use three special notations for integrals, depending on the role
of the measure. Expectations on the underlying probability space are
written $\ee{\cdot}$ and $\ee{\cdot \giv \cdot}$. For a deterministic
probability measure $\rho$ on a path space, we write
$\ee^{\rho}{\cdot}$ and $\ee^{\rho}{\cdot \giv \cdot}$ for
expectations of functionals of the coordinate process under $\rho$,
and use the same notation on the Polish space $S$ of
Lemma~\ref{l-pass}. For a random probability measure $M$
on $S$ and a bounded or nonnegative measurable function $F$ on $S$,
$\inner{M}{F}$ is the random variable
$\om \mapsto \int F \, \diffd M(\om)$. All other integrals are written with
$\int$. These include integrals in time and space, against measures
on Euclidean spaces, against random counting measures such as jump
measures, against mean measures, and against measures on $\sP(S)$.

Unless stated otherwise, filtrations satisfy the usual conditions.
Identities between processes hold up to indistinguishability. For stopping times $S$
and $T$, the sets $\llbracket S, T \rrbracket$,
$\llbracket S, T \llbracket$, $\rrbracket S, T \rrbracket$ and
$\rrbracket S, T \llbracket$ are the stochastic intervals in
$\Om \times [0,\infty)$. A process on $\llbracket 0, T \llbracket$ is
a family $\pseq*[ti]{Y_t}$ with $Y_t$ defined on $\set*{t < T}$. It is
adapted if each $Y_t \inds*{t < T}$ is $\sF_t$-measurable. We write
$Y_{t-}$ for left limits and $Y^T := \pseq*[ti]{Y_{t \wedge T}}$ for
the process stopped at $T$.

Compensators of random measures, and Poisson random measures
relative to a filtration, are defined as in
\cite[Section II.1]{JacShi03}. For a measurable space $U$,
predictable functions on $\Om \times (0,\infty) \times U$ are those
measurable for the product of the predictable $\ssig$ and the
$\ssig$ of $U$. Local times of semimartingales are symmetric local
times, normalized through Tanaka's formula with $\sgn(0) = 0$; see
\cite[Section IV.7]{Pro04}. For a
standard Brownian motion, the local time at $a$ is then the value at
$a$ of the continuous density of the occupation measure with respect
to Lebesgue measure. We refer to these as semimartingale local times
to distinguish them from the local times of
Section~\ref{ss-loctime}, which are normalized as occupation
densities.

\subsection{Use of generative AI}
During the preparation of this manuscript, the authors used
OpenAI's ChatGPT (GPT-5 and GPT-6) and Anthropic's Claude (Opus
5) to
assist with literature searches, to draft and polish exposition,
and to check the manuscript for errors.

Mathematically, such models have assisted us in adapting existing
proofs to somewhat weaker integrability assumptions on the
initial conditions.
In no case
did the main idea originate from AI. While the main cancellations
involved in the proof of Proposition~\ref{p-nonexpl} were
designed by the authors, AI suggested the exact exponential form
for the function $\psi$.

All AI-generated output, including the references identified
through literature searches, was independently reviewed and
verified by the authors. The authors take full responsibility for
the correctness, originality, and content of the manuscript,
including any errors that may remain.

\section{The model and its first properties}\label{s-model}

\subsection{Model and main result}\label{ss-main}
Particle $i \in [N]$ has a position $X^i$ and a jump intensity
$\ld^i \geq 0$, which evolve as
\begin{align*}
  \diffd X^i_t &= \diffd B^i_t + \int_E z \, \mu^i(\dt,\dz), \\
  \diffd\ld^i_t &= -\ld^i_t \, \dt + \gm_N \sum_{j \neq i} \diffd\ell^{i,j}_t .
\end{align*}
Here $B^1, \dots, B^N$ are independent Brownian motions, and $\mu^i$
is the jump measure of particle $i$. It has compensator
$\ld^i_t \, \dt \, \nu(\dz)$, so particle $i$ jumps at rate $\ld^i_t$,
by an amount with law $\nu$, a probability measure on
$E = \R \setminus \set{0}$. The collision local time
\begin{align*}
  \ell^{i,j}_t = \lim_{\eps \downarrow 0} \frac{1}{2\eps} \int_0^t
  \inds{\abs*{X^i_s - X^j_s} < \eps} \, \ds
\end{align*}
models the time that particles $i$ and $j$ have spent together, and
$\gm_N \geq 0$ is the interaction strength. The precise definition is
Definition~\ref{d-sol}. 

We are particularly interested in the mean-field limit of the
system, i.e., when $N \gm_N \to \kp \geq 0$. We
expect 
$\gm_N \sum_{j \neq i} \diffd\ell^{i,j}_t$ to approach
$\kp \, u(t, X^i_t) \, \dt$, where $u(t,\cdot)$ is the density of the
position of a typical particle at time $t$.
(Remark~\ref{r-meanfield}).
This leads to the \define{nonlinear system}
\begin{align}
  \label{n-sys-X}
  \diffd X_t &= \diffd B_t + \int_E z \, \mu(\dt,\dz), \\
  \label{n-sys-ld}
  \diffd\ld_t &= \prn{\kp \, u(t, X_t) - \ld_t} \, \dt ,
  \ewhere u(t,\cdot) = \text{ density of the law of } X_t,
\end{align}
where $\mu$ has compensator $\ld_t \, \dt \, \nu(\dz)$. Its solutions
are laws $P$ of the triple $(X, B, \ld)$ on the path space $\Om_1$ of
Section~\ref{ss-paths} (Definition~\ref{d-nl}).
For $N \geq 2$, let $\eta_N \in \sP\prn*{(\R \times [0,\infty))^N}$
be the joint law of the initial conditions
$(X^{i,N}_0, \ld^{i,N}_0)_{i \in [N]}$ of the $N$-particle system
with interaction strength $\gm_N$. 
With the empirical measure $Q_N$ defined by
\begin{align*}
  Q_N := \oo{N} \sum_{i=1}^N \delta_{(X^{i,N}, B^{i,N},
  \ld^{i,N})},
\end{align*}
we state our main result:

\begin{theorem}\label{t-intro}
  For every $N \geq 2$, $\gm_N \geq 0$ and $\eta_N$, the
  $N$-particle system has a solution, unique in law, which does not
  explode. 

  Suppose  that $N \gm_N \to \kp \geq 0$,
  that $\eta_N$ are exchangeable and $\eta_0$-chaotic for
  some $\eta_0 \in \sP\prn*{\R \times [0,\infty)}$, i.e., the law of
  $(X^{i,N}_0, \ld^{i,N}_0)_{i \in [k]}$ converges weakly to
  $\eta_0^{\otimes k}$ for each $k \in \N$, and that
  the random variables $\ld^{1,N}_0$, $N \geq 2$, are uniformly
  integrable. Then
  \begin{enumerate}[(a)]
  \item\label{i-intro-nl}
    the nonlinear
  system \eqref{n-sys-X}, \eqref{n-sys-ld} with the initial law $\eta_0$ has exactly one solution $P$, in the sense of Definition~\ref{d-nl};
  \item\label{i-intro-chaos}
    as $N \to \infty$, we have $Q_N \to P$
    in probability, and the law of
    $(X^{i,N}, B^{i,N}, \ld^{i,N})_{i \in [k]}$ converges weakly to
    $P^{\otimes k}$ for each $k\in\N$.
  \end{enumerate}
\end{theorem}

Theorem~\ref{t-main} also gives convergence of the pair empirical
measures, which include the collision local times of pairs of
particles, to the law of two independent copies of the limit together
with their collision local time. 

\subsection{Heuristics}\label{ss-heuristic}
The interaction is supported on pairwise encounters, but a particle potentially
meets many other particles. Under mean-field scaling, its
accumulated input should therefore average over their trajectories.
If particles become asymptotically independent with position
density $u$, this suggests
\begin{align*}
  \frac{1}{N-1} \sum_{j \neq i} \ell^{i,j}_t
  \approx \int_0^t u(s, X^i_s) \, \ds,
\end{align*}
as in the averaging mechanism underlying the Burgers model
\cite{Szn86}. In the present setting, the averaged input enters
the intensity and retains the exponential memory of past
encounters. Lemma~\ref{l-condlt} makes the conditional-local-time
identity precise, and Section~\ref{s-ident} justifies its use in
the particle limit.

Two features require estimates beyond a formal averaging argument.
First, the singular excitation could in principle lead to an
accumulation of jumps. The non-explosion proof uses a function of
pairwise distances together with the total intensity, chosen so
that the local-time contributions cancel. A related cancellation
gives a superlinear moment bound on individual intensities that is
uniform in the number of particles. Second, collision local time
is not directly a continuous functional of the paths in the
topology used for tightness. We therefore retain it in the pair
empirical measure and identify its limit through occupation-density
approximations and uniform integrability.

The limiting intensity still depends on the particle's own path,
so propagation of chaos does not eliminate the jumps. It replaces
the random input from other particles by a density-dependent
input. The density has a heat-kernel representation, which gives
finite-horizon bounds and permits a coupling of two candidate
solutions with the same Brownian motion and Poisson noise.
Gronwall's lemma then yields uniqueness of the nonlinear law.

\subsection{Piecewise Brownian processes}\label{ss-loctime}
Collision local times will be local times of
differences of particle positions. Between jumps, positions and
their differences are Brownian motions with a constant variance,
plus constants. We set up the
definitions in this generality, since they are also used for the
nonlinear system in Section~\ref{s-nonlin}.

\begin{definition}\label{d-pwb}
Let $c > 0$. A \cd{} process $Z$ on
a filtered probability space
$(\Om, \sF, \bF, \bP)$ which satisfies the usual conditions is
said to be \define{piecewise Brownian with variance $c$} if
\begin{align*}
Z = Z_0 + W + V,
\end{align*}
where
\begin{enumerate}[(i)]
\item\label{i-pwb-init} $Z_0$ is $\sF_0$-measurable,
\item\label{i-pwb-bm} $W / \sqrt{c}$ is a standard
  $\bF$-Brownian motion,
\item\label{i-pwb-jump} $V$ is a \cd{} $\bF$-adapted piecewise constant process
  with $V_0=0$ and finitely many jumps on bounded intervals.
\end{enumerate}
\end{definition}

For a piecewise Brownian process $Z = Z_0 + W + V$ with
variance $c$ and
$a \in \R$, we define
\begingroup
\medmuskip=3mu
\thickmuskip=4mu
\begin{equation}\label{tanaka}
  c \, \ell_t(a) := \abs*{Z_t - a} - \abs*{Z_0 - a}
  - \int_0^t \sgn(Z_{s-} - a) \, \diffd W_s
  - \sum_{0 < s \leq t} \prn{\abs*{Z_s - a} - \abs*{Z_{s-} - a}},
\end{equation}
\endgroup
and $\ell := \ell(0)$. Between jumps, $Z$ is a Brownian motion
with variance $c$
plus a constant. Hence it spends zero Lebesgue time at each level,
the value of $\sgn$ at $0$ in \eqref{tanaka} is immaterial, and
the right side of \eqref{tanaka} is the semimartingale local time of $Z$ at
$a$ \cite[Section IV.7]{Pro04}. So $\ell_t(a)$ is $1/c$ times
the semimartingale local
time of $Z$ at $a$. It is continuous and nondecreasing in $t$, starts
at $0$, increases only on $\set*{Z = a}$, and is jointly continuous in
$(t,a)$. Since $\qv{W}_t = ct$, it is the value at $a$ of the
continuous density of the occupation measure
$A \mapsto \int_0^t \ind{A}(Z_s) \, \ds$. Equivalently,
\begin{align*}
  \ell_t(a) = \lim_{\eps \downarrow 0} \frac{1}{2\eps}
  \int_0^t \inds{\abs*{Z_s - a} < \eps} \, \ds.
\end{align*}
We call $\ell_t(a)$ the \define{local time} of $Z$ at $a$. In
particular, $\ell_t(a)$ is a measurable function of the path of $Z$
on $[0,t]$, and the local time of $-Z$ at $a$ is that of $Z$
at $-a$. This normalization, as an occupation density, differs
from the semimartingale one unless $c = 1$.

Let $\zeta$ be a stopping time. A process $Z$ on
$\llbracket 0, \zeta \llbracket$ is \define{piecewise Brownian
with variance $c$
on $\llbracket 0, \zeta \llbracket$} if there are stopping times
$\zeta_n \uparrow \zeta$ and piecewise Brownian processes
$Z^{(n)}$ with variance $c$ with $Z = Z^{(n)}$ on
$\llbracket 0, \zeta_n \llbracket$.
Since local times depend only on the past, those of
$Z^{(n)}$ and $Z^{(m)}$ agree on
$\llbracket 0, \zeta_n \wedge \zeta_m \llbracket$. Together they
define the local time $\ell_t(a)$ of $Z$ for $t < \zeta$.

\subsection{Particle dynamics}\label{ss-dynamics}
Fix $N \in \N$ with $N \geq 2$, an interaction strength $\gm \geq 0$,
and a probability measure $\nu$ on the space of jump sizes
\begin{align*}
  E := \R \setminus \set{0}.
\end{align*}
We consider $N$ particles on the real line. Particle $i \in [N]$ has
a position $X^i$ and a jump intensity $\ld^i \geq 0$. Positions move
as independent Brownian motions, and particle $i$ jumps at rate
$\ld^i$ by an amount with law $\nu$. The intensity $\ld^i$ decays
exponentially and is pushed up whenever particle $i$ meets another
particle. Meetings are measured by collision local times, which
record the time two particles spend together. Since $N$ is fixed in
this
section, the dependence on $N$ is not recorded. When it matters, we
write $X^{i,N}$, $\ld^{i,N}$, and so on.

For an $N$-dimensional Brownian motion $B = (B^1, \dots, B^N)$ and
positions $X = (X^1, \dots, X^N)$, we set, for $i \neq j$,
\begin{equation}\label{pairs}
  Z^{i,j} := X^i - X^j
  \quad\text{and}\quad
  W^{i,j} := B^i - B^j.
\end{equation}
Each $W^{i,j}$ is $\sqrt{2}$ times a standard Brownian motion,
so
$\qv{W^{i,j}}_t = 2t$. Note
that $Z^{j,i} = -Z^{i,j}$ and $W^{j,i} = -W^{i,j}$.

\begin{definition}\label{d-sol}
  Let $\eta \in \sP\prn*{(\R \times [0,\infty))^N}$. A
  \define{solution up to explosion} of the $N$-particle system with
  interaction strength $\gm$, jump law $\nu$, and initial law $\eta$
  consists of
  \begin{enumerate}[(i)]
  \item\label{i-space} a filtered probability space $(\Om, \sF, \bF, \bP)$, with
    $\bF = \filt*$ satisfying the usual conditions, and an
    $N$-dimensional standard $\bF$-Brownian motion
    $B = (B^1, \dots, B^N)$;
  \item\label{i-init} $\sF_0$-measurable random variables
    $(X^i_0, \ld^i_0)_{i \in [N]}$ with joint law $\eta$;
  \item\label{i-mu} random measures $\mu^1, \dots, \mu^N$ on
    $(0,\infty) \times E$, each a sum of unit point masses with at
    most one atom per time, such that $t \mapsto \mu^i((0,t] \times A)$
    is $\bF$-adapted for every Borel $A \subseteq E$;
  \item\label{i-proc} $\bF$-adapted processes $X = (X^i)_{i \in [N]}$ and
    $\ld = (\ld^i)_{i \in [N]}$ on $\llbracket 0, \zeta \llbracket$,
    where
    \begin{equation}\label{Tn}
      J_t := \sum_{i=1}^N \mu^i\prn{(0,t] \times E}, \quad
      T_n := \inf\sets{t \geq 0}{J_t \geq n}, \quad
      \zeta := \sup_{n \in \N} T_n,
    \end{equation}
  \end{enumerate}
  such that the following conditions hold.
  \begin{enumerate}[(a)]
  \item\label{i-nosim} Almost surely, $\sum_{i=1}^N \mu^i\prn*{\set*{t} \times E} \leq 1$
    for all $t > 0$, and $\mu^i\prn*{[\zeta,\infty) \times E} = 0$ for
    all $i \in [N]$.
  \item\label{i-eqs} Almost surely, for all $i \in [N]$ and $t \in [0,\zeta)$,
    \begin{align}
      X^i_t &= X^i_0 + B^i_t + \int_{(0,t] \times E} z \, \mu^i(\ds,\dz),
      \label{X} \\
      \ld^i_t &= e^{-t} \ld^i_0
      + \gm \sum_{j \neq i} \int_0^t e^{-(t-s)} \, \diffd\ell^{i,j}_s,
      \label{ld}
    \end{align}
    where $\ell^{i,j}$ is the local time at $0$ of
    $Z^{i,j}$ on $\llbracket 0, \zeta \llbracket$, in the sense of
  Section~\ref{ss-loctime}, with $c = 2$ and $W = W^{i,j}$.
  \item\label{i-comp} For all $i \in [N]$ and $n \in \N$, the $\bF$-compensator of
    the random measure $\inds*{s \leq T_n} \mu^i(\ds,\dz)$ is
    \begin{align*}
      \inds{s \leq T_n} \ld^i_s \, \ds \, \nu(\dz).
    \end{align*}
  \end{enumerate}
  We call $\zeta$ the \define{explosion time} and $T_n$ the $n$-th
  \define{jump time}. The solution is \define{global} if
  $\zeta = \infty$ a.s.
\end{definition}

\begin{remark}\label{r-sol}
  \begin{enumerate}[(1)]
  \item\label{i-loctime} For $t < \zeta$ we have $J_t < \infty$,
    so the integral in \eqref{X} is a finite sum, and
    $Z^{i,j} - W^{i,j}$ is constant on each
    $\llbracket T_n, T_{n+1} \llbracket$. So $Z^{i,j}$ is piecewise
    Brownian with variance $2$ on $\llbracket 0, \zeta \llbracket$, with
    $\zeta_n := T_n$ and
    $Z^{(n)} := Z^{i,j}_0 + W^{i,j}
    + (Z^{i,j} - Z^{i,j}_0 - W^{i,j})^{T_n}$, and $\ell^{i,j}$ is well
    defined. If the solution is global, $Z^{i,j}$ is piecewise
    Brownian with variance $2$. We call $\ell^{i,j}$ the \define{collision local time}
    of particles $i$ and $j$. It is continuous and nondecreasing,
    starts at $0$, and increases only on
    $\sets*{t}{Z^{i,j}_t = 0}$.
  \item\label{i-occup} By Section~\ref{ss-loctime}, $\ell^{i,j}$
    is half the semimartingale local time of $Z^{i,j}$ at $0$, $\ell^{i,j} =
    \ell^{j,i}$, and
    \begin{align*}
      \ell^{i,j}_t = \lim_{\eps \downarrow 0} \frac{1}{2\eps}
      \int_0^t \inds{\abs*{X^i_s - X^j_s} < \eps} \, \ds.
    \end{align*}
    In particular, $\ell^{i,j}$ is a measurable function of the path
    of $X^i - X^j$.
  \item\label{i-lddyn} Equation \eqref{ld} is equivalent to
    $\diffd\ld^i_t = -\ld^i_t \, \dt + \gm \sum_{j \neq i} \diffd\ell^{i,j}_t$. Each
    $\ld^i$ is continuous and nonnegative, and
    $t \mapsto e^t \ld^i_t$ is nondecreasing on $[0,\zeta)$.
  \item\label{i-rate} Condition~\ref{i-comp} says that particle $i$ jumps at rate $\ld^i_t$
    and that its jump sizes have law $\nu$. By \ref{i-nosim}, two particles
    never jump at the same time.
  \item\label{i-scaling} The decay rate of $\ld^i$ and the diffusion coefficient are
    both $1$. This is no loss of generality. If $\ld^i$ decays at
    rate $a > 0$, the change of variables $t' = a t$ and
    $X'_{t'} = \sqrt{a}\, X_{t'/a}$ leads to the present system,
    with $X_0$ replaced by $\sqrt{a}\, X_0$, $\ld_0$ by $\ld_0/a$,
    $\nu$ by its image under $z \mapsto \sqrt{a}\, z$, and $\gm$ by
    $\gm a^{-3/2}$. A diffusion
    coefficient $\sg \neq 1$ is removed in the same way by
    $X \mapsto X/\sg$, with $X_0$ replaced by $X_0/\sg$, $\nu$ by its
    image under $z \mapsto z/\sg$, and $\gm$ by $\sg \gm$, when
    $\ell^{i,j}$ is still defined by \eqref{tanaka} with $c = 2$. Since
    $\ell^{i,j}$ is half the semimartingale local time of $Z^{i,j}$, it scales like
    $Z^{i,j}$.
  \end{enumerate}
\end{remark}

\subsection{Existence up to explosion}\label{ss-exist}
We construct a solution by Poisson thinning.
Let
$(\Om, \sF, \bP)$ be a probability space that supports the following
independent objects:
\begin{enumerate}[(1)]
\item\label{i-drv-init} a random vector $(X^i_0, \ld^i_0)_{i \in [N]}$ with law
  $\eta \in \sP\prn*{(\R \times [0,\infty))^N}$;
\item\label{i-drv-bm} an $N$-dimensional standard Brownian motion
  $B = (B^1, \dots, B^N)$;
\item\label{i-drv-prm} Poisson random measures $\Pi^1, \dots, \Pi^N$ on
  $(0,\infty) \times E \times (0,\infty)$, each with intensity
  $\ds \, \nu(\dz) \, \dy$.
\end{enumerate}
Let $\bF = \filt*$ be the usual augmentation of the filtration
\begin{align*}
  \sF^0_t := \sig{X^i_0, \ld^i_0, B_s, \Pi^i(A) : i \in [N],\
  s \leq t,\ A \subseteq (0,t] \times E \times (0,\infty)
  \text{ Borel}}.
\end{align*}
It is standard that $B$ is an $\bF$-Brownian motion and that each
$\Pi^i$ is an $\bF$-Poisson random measure. In particular, the
$\bF$-compensator of $\Pi^i$ is $\ds \, \nu(\dz) \, \dy$
\cite[Section II.1]{JacShi03}.

\begin{proposition}\label{p-exist}
  In the setting above, there exist processes $X$, $\ld$ and random
  measures $\mu = (\mu^1, \dots, \mu^N)$ such that
  $(\Om, \sF, \bF, \bP; B, X, \ld, \mu)$ is a solution up to
  explosion with initial values $(X^i_0, \ld^i_0)_{i \in [N]}$, in
  the sense of Definition~\ref{d-sol}. Moreover, the following hold
  almost surely.
  \begin{enumerate}[(a)]
  \item\label{i-thin} For each $i \in [N]$, with the integrand set to $0$ for
    $s \geq \zeta$,
    \begin{equation}\label{thin}
      \mu^i(\ds,\dz) = \int_{(0,\infty)} \inds{s < \zeta,\ y \leq \ld^i_s}
      \, \Pi^i(\ds,\dz,\dy).
    \end{equation}
  \item\label{i-blowup} On $\set*{\zeta < \infty}$, we have
    $\sup_{t < \zeta} \sum_{i=1}^N \ld^i_t = \infty$.
  \end{enumerate}
\end{proposition}

The proof is given in Appendix~\ref{ap-exist}.

\begin{remark}\label{r-thin-unique}
  The construction is pathwise, and the solution of
  Proposition~\ref{p-exist} is the only one driven by
  $\Pi^1, \dots, \Pi^N$. More precisely, let $(\tX, \tilde\ld,
  \tilde\mu)$ be another solution up to explosion on the same space,
  with the same Brownian motion and initial values, which satisfies
  \eqref{thin} with the same $\Pi^i$ and with
  $(\tilde\ld, \tilde\zeta)$ in place of $(\ld, \zeta)$. Then it has
  the same explosion
  time and coincides with $(X, \ld, \mu)$ on
  $\llbracket 0, \zeta \llbracket$. This follows by induction over the
  jump times. As long as the two solutions agree, so do their
  intensities, and \eqref{thin} forces both to accept the same next
  atom.
\end{remark}

\subsection{Non-explosion}\label{ss-nonexpl}
Let $\Ld := \sum_{i=1}^N \ld^i$ be the total intensity. By
condition~\ref{i-comp} of Definition~\ref{d-sol} and $\nu(E) = 1$, the
compensator of $J$ on $\llbracket 0, T_n \rrbracket$ is
$\int_0^\cdot \Ld_s \, \ds$. By \eqref{ld},
\begin{equation}\label{Ld}
  \diffd\Ld_t = -\Ld_t \, \dt + \gm \sum_{i \neq j} \diffd\ell^{i,j}_t.
\end{equation}

\begin{proposition}\label{p-nonexpl}
  Every solution up to explosion is global.
\end{proposition}

\begin{proof}
  Let $a := \gm N$, $\psi(w) := e^{-a\abs*{w}}$, and
  \begin{align*}
    V_t := \ot \sum_{i \neq j} \psi(Z^{i,j}_t) + N \Ld_t,
    \quad t < \zeta.
  \end{align*}
  Fix $n \in \N$ and work on $\llbracket 0, T_n \rrbracket$. The
  function $\psi$ is a difference of convex functions, with
  $\psi'' = a^2 \psi - 2a \delta_0$ in the sense of distributions.
  Since $\qv{W^{i,j}}_t = 2t$ and $\ell^{i,j}$ is half the semimartingale local time of
  $Z^{i,j}$ at $0$, the Meyer--\ito\ formula
  \cite[Section IV.7]{Pro04} gives
  \begin{multline*}
    \psi(Z^{i,j}_t) = \psi(Z^{i,j}_0)
    + \int_0^t \psi'(Z^{i,j}_{s-}) \, \diffd W^{i,j}_s
    + a^2 \int_0^t \psi(Z^{i,j}_s) \, \ds \\
    + \sum_{0 < s \leq t} \prn{\psi(Z^{i,j}_s) - \psi(Z^{i,j}_{s-})}
    - 2a \, \ell^{i,j}_t.
  \end{multline*}
  At each jump time only the $2(N-1)$ ordered pairs containing the
  jumping particle move, and each term changes by less than $1$, since
  $0 < \psi \leq 1$. We combine this with \eqref{Ld} and note that
  our choice of the constant $a$ makes the local-time
  terms cancel. Hence
  \begin{align*}
    V_{t \wedge T_n} \leq V_0 + M_{t \wedge T_n}
    + a^2 \binom{N}{2} \, t - \int_0^{t \wedge T_n} \Ld_s \, \ds,
  \end{align*}
  where
  \begin{align*}
    M := \ot \sum_{i \neq j} \int_0^\cdot \psi'(Z^{i,j}_{s-})
    \, \diffd W^{i,j}_s
    + (N-1) \prn{J - \int_0^\cdot \Ld_s \, \ds}
  \end{align*}
  is a martingale on $\llbracket 0, T_n \rrbracket$, because
  $\abs*{\psi'} \leq a$ and $J_{T_n} \leq n$. Since $V \geq 0$ and
  $V_0 \leq \binom{N}{2} + N \Ld_0$, taking conditional expectations
  gives
  \begin{align*}
    \ee{J_{t \wedge T_n} \giv \sF_0}
    = \ee{\int_0^{t \wedge T_n} \Ld_s \, \ds \giv \sF_0}
    \leq c_t := \binom{N}{2} \prn{1 + a^2 t} + N \Ld_0 .
  \end{align*}
  As $J_{t \wedge T_n} = n$ on
  $\set*{T_n \leq t}$, it follows that
  \begin{align*}
    \pp{\zeta \leq t \giv \sF_0} \leq \pp{T_n \leq t \giv \sF_0}
    \leq \frac{c_t}{n}.
  \end{align*}
  Since $n$ is arbitrary, $\zeta = \infty$ a.s.
\end{proof}

\subsection{Uniqueness in law}\label{ss-unique}
Remark~\ref{r-thin-unique} compares solutions driven by the same
Poisson random measures. Every solution is driven by such measures
on an extension of its probability space, and this gives uniqueness
in law.

\begin{proposition}\label{p-unique}
  Let $(\Om, \sF, \bF, \bP; B, X, \ld, \mu)$ be a solution up to
  explosion with initial law $\eta$. On an extension of
  $(\Om, \sF, \bP)$ there exist Poisson random measures
  $\Pi^1, \dots, \Pi^N$ on $(0,\infty) \times E \times (0,\infty)$,
  each with intensity $\ds \, \nu(\dz) \, \dy$, such that
  $(X_0, \ld_0)$, $B$, $\Pi^1, \dots, \Pi^N$ are independent and
  \eqref{thin} holds. In particular, the law of
  $(X_0, \ld_0, B, X, \ld, \mu)$ depends only on $\gm$, $\nu$ and
  $\eta$.
\end{proposition}

The proof is given in Appendix~\ref{ap-unique}.

\section{Uniform bounds and tightness}\label{s-bounds}\label{s-tight}

From now on the number of particles varies, and we work under
Assumption~\ref{a-std} below. It imposes the mean-field scaling
$N \gm_N \to \kp$ on the interaction strength. This section proves
bounds, uniform in $N$, on finite time horizons: a superlinear moment
bound on the running maximum of the intensities, uniform integrability
of the jump counts and their compensators, an Aldous-type
bound on the increments of the intensities, the positions and the
collision local times, and uniform integrability of the collision
local times. These estimates then yield tightness of the particle
and pair empirical measures, together with the product structure
of their limit points.

\subsection{Setting}\label{ss-bset}
For each $N \geq 2$ we are given an interaction strength $\gm_N$ and
an initial law $\eta_N \in \sP\prn*{(\R \times [0,\infty))^N}$. The
jump law $\nu$ does not depend on $N$. We write
$(x_i, y_i)_{i \in [N]}$ for the generic point of
$(\R \times [0,\infty))^N$. The following standing assumption on
these ingredients is in force for the rest of the paper.

\begin{assumption}\label{a-std}
  \begin{enumerate}[(i)]
  \item\label{i-mf} $N \gm_N \to \kp$ as $N \to \infty$, for some
    $\kp \geq 0$.
  \item\label{i-eta} Each $\eta_N$ is exchangeable, and
    the coordinate $y_1$ is uniformly integrable under
    $\seq[N]{\eta_N}$:
    \begin{equation}\label{ld0ui}
      \sup_{N \geq 2} \int (y_1 - m)^+ \, \diffd\eta_N \to 0
      \quad \text{as } m \to \infty .
    \end{equation}
  \item\label{i-chaos} There is $\eta_0 \in \sP(\R \times [0,\infty))$
    such that $\seq[N]{\eta_N}$ is $\eta_0$-chaotic: for each $k \in \N$, the
    law of $(x_i, y_i)_{i \in [k]}$ under $\eta_N$ converges weakly to
    $\eta_0^{\otimes k}$ as $N \to \infty$.
  \end{enumerate}
\end{assumption}

Given \ref{i-chaos}, condition \eqref{ld0ui} is equivalent to
$\int y_1 \, \diffd\eta_N \to \int y \, \diffd\eta_0 < \infty$.
The construction in Appendix~\ref{a-ui-function}, a smooth concave
version of de la Vall\'ee-Poussin's criterion, gives a function
\begin{align*}
  G(x) := \int_0^x g(u) \, \diffd u, \quad x \geq 0,
\end{align*}
where $g \colon [0,\infty) \to [0,\infty)$ is nondecreasing with
$g(0) = 0$ and $g(x) \to \infty$
as $x \to \infty$, such that
\begin{equation}\label{Gbound}
  \sup_{N \geq 2} \int G(y_1) \, \diffd\eta_N \leq 1 ,
\end{equation}
with $g$ concave and continuously differentiable on $[0,\infty)$.
We fix such a $G$ for the rest of the paper. It is convex,
continuous and nondecreasing, with $G(0) = 0$, and
$G(x)/x \geq g(x/2)/2 \to \infty$ as $x \to \infty$. Concavity
of $g$ and $g(0) = 0$ give $g(u) \geq (u/x) g(x)$ for
$0 \leq u \leq x$, $g(2x) \leq 2 g(x)$, and $x g'(x) \leq g(x)$.
Integration and convexity of $G$ then give, for $x > 0$ and $c \geq 1$,
\begin{equation}\label{Gprop}
  x g'(x) \leq g(x) \leq \frac{2 G(x)}{x}, \qquad
  G(x) \leq x g(x), \qquad
  G(cx) \leq c^2 G(x) .
\end{equation}
Indeed, $g(cu) \leq c g(u)$ for $c \geq 1$ by concavity and
$g(0)=0$, so
\begin{align*}
  G(cx) = c \int_0^x g(cu) \, \diffd u
  \leq c^2 \int_0^x g(u) \, \diffd u = c^2 G(x).
\end{align*}
Moreover, for $x, y \geq 0$, the concave function
$w \mapsto g(x) w - G(w)$ is maximal at $w = x$, so
\begin{equation}\label{Gyoung}
  g(x) \, y \leq G(y) + x g(x) - G(x) \leq G(y) + x g(x) .
\end{equation}
Finally, since $G(x)/x \to \infty$, there is $x_G \geq 0$ with
$x \leq x_G + G(x)$ for all $x \geq 0$, and a family of nonnegative
random variables $\xi$ with $\sup \ee{G(\xi)} < \infty$ is uniformly
integrable, because
$\ee{\xi \inds{\xi \geq m}} \leq \sup_{x \geq m} (x/G(x)) \, \ee{G(\xi)}$.

For $N \geq 2$, let
$(X^N, \ld^N, \mu^N)$ be a solution with interaction strength $\gm_N$,
jump law $\nu$ and initial law $\eta_N$. By
Propositions~\ref{p-nonexpl} and~\ref{p-unique}, it is global and its
law is determined by $\gm_N$, $\nu$ and $\eta_N$. We write $X^{i,N}$,
$B^{i,N}$, $\ld^{i,N}$, $\mu^{i,N}$, $\ell^{i,j,N}$, $J^N$ and
$\Ld^N := \sum_{i=1}^N \ld^{i,N}$, and we drop the index $N$ when no
confusion arises. Each $\ld^{i,N}$ is driven by $N - 1$ collision
local times with coefficient $\gm_N$, so the total interaction
strength felt by one particle is
$(N-1) \gm_N \leq \bar\kp$ where
\begin{align*}
  \bar\kp := \sup_{N \geq 2} N \gm_N < \infty.
\end{align*}
Relabeling the particles preserves the equations and the
exchangeable initial law. Uniqueness in law
(Proposition~\ref{p-unique}) therefore implies that the family
$(X^{i,N}, B^{i,N}, \ld^{i,N}, \mu^{i,N})_{i \in [N]}$
is exchangeable. By
Definition~\ref{d-sol}\ref{i-init} and exchangeability,
$\ee*{G(\ld^{i,N}_0)} = \int G(y_1) \, \diffd\eta_N \leq 1$ for all $i \in [N]$. By
\ref{i-chaos}, the continuity of $G$ and Fatou's lemma, also
\begin{equation}\label{G0}
  \int G(y) \, \diffd\eta_0 \leq 1,
  \eand \int y \, \diffd\eta_0 \leq x_G + 1 .
\end{equation}

For $T \geq 0$, $C_T$ denotes a constant that depends only on $T$,
$\bar\kp$ and the function $G$. In
particular, $C_T$ does not depend on $N$.
It can always be chosen to be nondecreasing in $T$, 
and its value may change from use to use.

By Remark~\ref{r-sol}\ref{i-occup}, $\ell^{i,j,N}$ is a measurable
function of $X^{i,N} - X^{j,N}$. Hence the law of
$(X^{i,N}, \allowbreak B^{i,N}, \allowbreak \ld^{i,N}, \allowbreak \mu^{i,N}, \allowbreak X^{j,N}, \allowbreak B^{j,N}, \allowbreak
\ld^{j,N}, \allowbreak \mu^{j,N}, \allowbreak \ell^{i,j,N})$ does not depend on the pair
$i \neq j$.

\subsection{Moments}\label{ss-bmom}

\begin{proposition}\label{p-moments}
  For all $N \geq 2$, $i \in [N]$ and $T \geq 0$,
  \begin{align*}
    \ee{G\prn{\textstyle\sup_{s \leq T} \ld^i_s}} \leq C_T .
  \end{align*}
\end{proposition}

\begin{proof}
  By Remark~\ref{r-sol}\ref{i-lddyn}, $s \mapsto e^s \ld^i_s$ is
  nondecreasing, so $\sup_{s \leq T} \ld^i_s \leq e^T \ld^i_T$
  and, by \eqref{Gprop}, it
  suffices to show that $\ee*{G(\ld^i_T)} \leq C_T$.

  Let $a := 2 (N-1) \gm$, so that $a/(2(N-1)) = \gm$
  and $a \leq 2\bar\kp$. With $\psi(w) := e^{-a\abs*{w}}$,
  define
  \begin{align*}
    \Psi := \oo{N-1} \sum_{j \neq i} \psi(Z^{i,j}) \in (0,1],
    \qquad
    \rho := \frac14 \, g(\ld^i) \, \Psi + G(\ld^i) .
  \end{align*}
  Let $\sg_n := \inf\sets*{t \geq 0}{\max_k \ld^k_t \geq n}$. It is a
  symmetric function of the particles, and $\sg_n \to \infty$ because
  the intensities are continuous and the solution is global.

  By the Meyer--\ito\ formula, as in the proof of
  Proposition~\ref{p-nonexpl}, the chain rule
  $\diffd g(\ld^i) = g'(\ld^i) \, \diffd\ld^i$ and
  $\diffd G(\ld^i) = g(\ld^i) \, \diffd\ld^i$ ($\ld^i$ is continuous
  and of finite variation, and $g \in C^1$), and the product rule
  ($g(\ld^i)$ is
  continuous and of finite variation),
  \begin{align*}
    \diffd\rho = \diffd M + \frac14 \, g(\ld^i) \Big( a^2 \Psi \, \dt
    - \frac{2a}{N-1} \tsum_{j \neq i} \diffd\ell^{i,j} \Big)
    + \frac14 \, \Psi \, g'(\ld^i) \, \diffd\ld^i
    + g(\ld^i) \, \diffd\ld^i + \diffd\Sigma,
  \end{align*}
  where $M$ is a stochastic integral against the $W^{i,j}$ and
  $\Sigma$ is the sum of the jumps of $\rho$.
  We define $\tau := T \wedge \sg_n$ and observe that
  the integrand of $M$ is
  bounded on $\rrbracket 0, \tau \rrbracket$, so $\ee*{M_\tau} = 0$.

  \emph{Drift and local times.} Since $\Psi \leq 1$, we have
  $a^2 g(\ld^i) \Psi \leq a^2 g(\ld^i)$. By \eqref{ld}, $g' \geq 0$
  and \eqref{Gprop},
  \begin{align*}
    \Psi \, g'(\ld^i) \, \diffd\ld^i
    \leq g'(\ld^i) \, \gm \sum_{k \neq i} \diffd\ell^{i,k}
    = g'(\ld^i) \prn{\diffd\ld^i + \ld^i \, \dt}
    \leq \diffd g(\ld^i) + g(\ld^i) \, \dt .
  \end{align*}
  Also
  $g(\ld^i) \, \diffd\ld^i = -\ld^i g(\ld^i) \, \dt
  + \gm \, g(\ld^i) \sum_{k \neq i} \diffd\ell^{i,k}$, and the last term cancels
  the local-time term in the display, since $a/(2(N-1)) = \gm$.

  \emph{Jumps.} Since $\ld$ is continuous, $\rho$ jumps only through
  $\Psi$. A jump of particle $i$ raises $\Psi$ by less than $1$, and
  a jump of particle $j \neq i$ by less than $1/(N-1)$. By
  condition~\ref{i-comp} of Definition~\ref{d-sol},
  \begin{align*}
    \ee{\Sigma_\tau}
    \leq \frac14 \, \ee{\int_0^\tau g(\ld^i_s) \Big( \ld^i_s
    + \oo{N-1} \sum_{j \neq i} \ld^j_s \Big) \ds}
    \leq \frac34 \, \ee{\int_0^\tau \ld^i_s \, g(\ld^i_s) \, \ds},
  \end{align*}
  where the second inequality uses
  $g(\ld^i) \ld^j \leq \ld^i g(\ld^i) + G(\ld^j)$ from \eqref{Gyoung},
  $G(\ld^j) \leq \ld^j g(\ld^j)$ and the fact that
  $\ee*{\int_0^\tau \ld^j_s \, g(\ld^j_s) \, \ds}$ does not depend on $j$, by
  exchangeability and the symmetry of $\sg_n$.

  \emph{Conclusion.} All terms below are finite. On
  $\set*{\tau>0}$, $\ld^k_s\leq n$ for $s\leq\tau$ and all $k$;
  on $\set*{\tau=0}$, the time integrals vanish and
  $\ld^i_\tau=\ld^i_0$. The endpoint terms are integrable because
  $\ee{G(\ld^i_0)}\leq1$ and
  $\ee{g(\ld^i_0)}\leq g(1)+2\ee{G(\ld^i_0)}\leq g(1)+2$ by
  \eqref{Gprop} and \eqref{Gbound}. Collecting the
  bounds, in which $\int_0^\tau \ld^i_s \, g(\ld^i_s) \, \ds$ has net
  coefficient $-1 + \frac34 = -\frac14$, we get
  \begin{align*}
    \ee{\rho_\tau} + \frac14 \, \ee{\int_0^\tau \ld^i_s \, g(\ld^i_s) \, \ds}
    \leq \ee{\rho_0} + \frac{a^2 + 1}{4} \, \ee{\int_0^\tau g(\ld^i_s) \, \ds}
    + \frac14 \, \ee{g(\ld^i_\tau)}.
  \end{align*}
  We absorb the last two terms into the left side. First,
  $b g(x) \leq x g(x) + b g(b)$ for $x\geq0$ and $b > 0$, so, with $b := a^2 + 1$,
  $\frac{a^2+1}{4} g(x) \leq \frac14 x g(x) + \frac{a^2+1}{4} g(a^2 + 1)$.
  Second,
  $g(x) \leq 2 G(x) + g(1)$ by \eqref{Gprop}, and $G(\ld^i) \leq \rho$.
  Hence
  \begin{align*}
    \ee{G(\ld^i_\tau)} \leq \ee{\rho_\tau}
    \leq 2 \, \ee{\rho_0} + \ot (a^2 + 1) \, g(a^2 + 1) \, T + \ot g(1) .
  \end{align*}
  The right side is at most $C_T$, since
  $\rho_0 \leq \frac14 g(\ld^i_0) + G(\ld^i_0)$, by \eqref{Gbound},
  \eqref{Gprop} and $a \leq 2\bar\kp$.
  Letting $n \to \infty$, Fatou's lemma and the continuity of $G$
  give
  $\ee*{G(\ld^i_T)} \leq C_T$.
\end{proof}

\begin{corollary}
\label{c-mom}
Let $T \geq 0$.
\begin{enumerate}[(a)]
\item\label{i-mom-sup} $\ee*{\sup_{s \leq T} \ld^i_s} \leq C_T$ for all
  $N \geq 2$ and $i \in [N]$, 
\item\label{i-mom-jump} $\ee*{J^N_T} = N \, \ee*{\int_0^T \ld^{1,N}_s \, \ds}
\leq N C_T$ for all $N \geq 2$, and 
\item\label{i-mom-ui} the family of all $\ld^{i,N}_\tau$, with $N \geq 2$,
$i \in [N]$ and $\tau \leq T$ a stopping time, satisfies
$\ee*{G(\ld^{i,N}_\tau)} \leq C_T$, and is therefore uniformly integrable.
\end{enumerate}
\end{corollary}

\begin{proof}
  Parts \ref{i-mom-sup} and \ref{i-mom-ui} follow from
  Proposition~\ref{p-moments}, since $x \leq x_G + G(x)$,
  $G$ is nondecreasing, and a family with bounded $G$-moments is
  uniformly integrable (Section~\ref{ss-bset}). For \ref{i-mom-jump}, the compensator
  of $J^N$ is $\int_0^\cdot \Ld^N_s \, \ds$, by
  condition~\ref{i-comp} of Definition~\ref{d-sol}, and
  $\ee*{\int_0^T \ld^{i,N}_s \, \ds}$ does not depend on $i$, by
  exchangeability.
\end{proof}

For $N \geq 2$, $i \in [N]$ and $t \geq 0$, let
$J^{i,N}_t := \mu^{i,N}((0,t] \times E)$ be the number of jumps of
$X^{i,N}$ in $(0,t]$.
Lemma~\ref{l-jumpmom} bounds the expected number of jumps with
marks in any Borel set and proves uniform integrability of the
jump counts and their compensators. Part~\ref{i-bdnu} makes the
probability of a jump larger than $a$ tend to zero, uniformly in $N$,
as $a\to\infty$. Small jumps satisfy
$\int_{(0,T]\times\set*{|z|<\dl}} |z|\,\mu^{1,N}(\ds,\dz)
\leq\dl J^{1,N}_T$; their expected total absolute size is at most
$\dl C_T$, by Corollary~\ref{c-mom}\ref{i-mom-jump} and exchangeability.

\begin{lemma}\label{l-jumpmom}
  Let $N \geq 2$ and $T \geq 0$.
  \begin{enumerate}[(a)]
    \item \label{i-bdnu} $\ee*{\mu^{1,N}((0,T] \times A)} \leq C_T \, \nu(A)$, 
    for every Borel $A \subseteq E$,
  \item \label{i-bdjsq} the families $\seq[N]{J^{1,N}_T}$ and
  $\seq[N]{\int_0^T \ld^{1,N}_s \, \ds}$ are uniformly integrable.
  \end{enumerate}
\end{lemma}
\begin{proof}
  By condition~\ref{i-comp} of Definition~\ref{d-sol}
  and Corollary~\ref{c-mom}\ref{i-mom-sup},
  we have
  \begin{align*}
    \ee*{\mu^{1,N}((0,T] \times A)}
    = \nu(A) \, \ee*{\int_0^T \ld^{1,N}_s \, \ds}
    \leq C_T \, \nu(A).
  \end{align*}
  For $T=0$, both families in \ref{i-bdjsq} consist of zero.
  For $T > 0$, Jensen's inequality gives
  \begin{align*}
  G\prn{\oo{T} \int_0^T \ld^{1,N}_s \, \ds}
  \leq \oo{T} \int_0^T G(\ld^{1,N}_s) \, \ds
  \leq G\prn{\sup_{s \leq T} \ld^{1,N}_s},
  \end{align*}
  whose right side has expectation at
  most $C_T$ by Proposition~\ref{p-moments}. So the second family in
  \ref{i-bdjsq} is uniformly integrable (Section~\ref{ss-bset}).
  The compensator
  $\int_0^\cdot \ld^{1,N}_s \, \ds$ of $J^{1,N}$ is continuous, so
  $\ee*{(J^{1,N}_T - \int_0^T \ld^{1,N}_s \, \ds)^2}
  = \ee*{\int_0^T \ld^{1,N}_s \, \ds}$$\leq C_T$, and
  $J^{1,N}_T$ is the sum of $\int_0^T \ld^{1,N}_s \, \ds$ and a
  random variable bounded in $\bL^2$, uniformly in $N$. This gives \ref{i-bdjsq}.
\end{proof}

\subsection{Increments}\label{ss-bincr}
The key local-time estimate is Lemma~\ref{l-pieces} in
Appendix~\ref{a-localtimes}. It bounds the expected increase in
local time over a stopping interval by the square root of its
length times the expected number of Brownian pieces in that
interval. Its other parts give uniform integrability and control
the change in local time when its spatial level varies. These
three conclusions supply, respectively, the increment estimate
below, Corollary~\ref{c-ltui}, and the local-time identification
in Proposition~\ref{p-ipair}.

\begin{proposition}\label{p-aldous}
  Let $N \geq 2$, $T \geq 0$ and $\theta \in [0,1]$, and let
  $\sg \leq \tau \leq \sg + \theta$ be stopping times with $\sg \leq T$. Then,
  for all $i \neq j$ in $[N]$,
  \begin{align*}
    \ee{1 \wedge \abs*{X^i_{\tau} - X^i_{\sg}}}
    + \ee{\ell^{i,j}_{\tau} - \ell^{i,j}_{\sg}}
    + \ee{\abs*{\ld^i_{\tau} - \ld^i_{\sg}}}
    \leq C_T \sqrt{\theta} .
  \end{align*}
\end{proposition}

\begin{proof}
  Since $\tau \leq T + 1$ and $\tau - \sg \leq \theta$,
  Corollary~\ref{c-mom}\ref{i-mom-sup} gives, for every $k \in [N]$,
  \begin{equation}\label{ldint}
    \ee{\int_{\sg}^{\tau} \ld^k_s \, \ds}
    \leq \theta \, \ee{\sup_{s \leq T+1} \ld^k_s}
    \leq C_T \theta .
  \end{equation}

  \emph{Positions.} 
  We have
\begin{align*}
  1 \wedge \abs*{X^i_{\tau} - X^i_{\sg}} \leq \abs*{B^i_{\tau} - B^i_{\sg}}
  + \mu^i((\sg,\tau] \times E).
\end{align*}
 By condition~\ref{i-comp} of
  Definition~\ref{d-sol} and \eqref{ldint}, the expectation of the
  right side is at most $\sqrt{\theta} + C_T \theta$.

  \emph{Local times.} By Proposition~\ref{p-nonexpl} and
  Remark~\ref{r-sol}\ref{i-loctime}, $Z^{i,j}$ is piecewise
  Brownian with variance $2$, with $W = W^{i,j}$ and $V = Z^{i,j} - Z^{i,j}_0 - W^{i,j}$. So
  Lemma~\ref{l-pieces}
  applies to $Z = Z^{i,j}$, with $\ell(0) = \ell^{i,j}$. The jumps of $V$ are those of $X^i$ and
  $X^j$, so the random variable $K$, defined in \eqref{num-pieces},
  satisfies
  $K \leq 1 + \mu^i((\sg,\tau] \times E) + \mu^j((\sg,\tau] \times E)$. By
  Lemma~\ref{l-pieces}\ref{i-pieces-bd}, condition~\ref{i-comp} of
  Definition~\ref{d-sol} and \eqref{ldint},
  \begin{align*}
    \ee{\ell^{i,j}_{\tau} - \ell^{i,j}_{\sg}}
    \leq \sqrt{\theta} \, \Big( 1 + \ee{\int_{\sg}^{\tau}
    \prn{\ld^i_s + \ld^j_s} \ds} \Big)
    \leq C_T \sqrt{\theta} .
  \end{align*}

  \emph{Intensities.} By \eqref{ld}, and since the $\ell^{i,k}$ are
  nondecreasing,
  \begin{align*}
    \abs*{\ld^i_{\tau} - \ld^i_{\sg}} \leq \int_{\sg}^{\tau} \ld^i_s \, \ds
    + \gm \sum_{k \neq i} \prn{\ell^{i,k}_{\tau} - \ell^{i,k}_{\sg}}.
  \end{align*}
  By \eqref{ldint}, the bound on the local times and
  $(N-1) \gm \leq \bar\kp$, the expectation of the right side is at
  most $C_T \theta + \bar\kp C_T \sqrt{\theta} \leq C_T \sqrt{\theta}$.
\end{proof}

\begin{corollary}\label{c-ltui}
  For every $T \geq 0$, the 
  sequence $\seq*[N]{\ell^{1,2,N}_T}$ is
  uniformly integrable.
\end{corollary}

\begin{proof}
  Apply Lemma~\ref{l-pieces}\ref{i-pieces-ui} to $Z = Z^{1,2}$, as in
  the proof of Proposition~\ref{p-aldous}, with $\sg = 0$,
$\tau = \theta = T$ and $\rho = \delta_0$. Here $K$, defined in
  \eqref{num-pieces}, satisfies
  $K \leq 1 + J^{1,N}_T + J^{2,N}_T$, so $\ee*{(K - M)^+} \leq \vp(M)
  := \sup_{N \geq 2} \ee*{(1 + J^{1,N}_T + J^{2,N}_T - M)^+}$, where
  $\vp(M) \to 0$ as $M \to \infty$ by
  Lemma~\ref{l-jumpmom}\ref{i-bdjsq} and exchangeability.
\end{proof}

\subsection{Path spaces}\label{ss-paths}
Let $C_\uparrow([0,\infty))$ be the set of continuous nondecreasing
$l \colon [0,\infty) \to [0,\infty)$ with $l_0 = 0$. It is a closed
subset of $C([0,\infty),\R)$. Define
\begin{align*}
  \Om_1 := D([0,\infty),\R) \times C([0,\infty),\R)
  \times C([0,\infty),[0,\infty)),
\end{align*}
and $\Om_2 := \Om_1 \times \Om_1 \times C_\uparrow([0,\infty))$
with Polish product topologies (\cite[Theorem
VI.1.14]{JacShi03});
we call them 
the \emph{canonical spaces} of \emph{one particle} and of a
\emph{pair}, respectively. 
For
$N \geq 2$ and
$i \neq j$ in $[N]$, let
\begin{align*}
  Y^{i,N} := (X^{i,N}, B^{i,N}, \ld^{i,N}) \in \Om_1,\ 
  \bar Y^{i,j,N} := (Y^{i,N}, Y^{j,N}, \ell^{i,j,N}) \in \Om_2, 
\end{align*}
where $\ell^{i,j,N} \in C_\uparrow([0,\infty))$ by
Remark~\ref{r-sol}\ref{i-loctime}. 
The \emph{empirical} and \emph{pair empirical}
measures are the random elements of the spaces
$\sP(\Om_1)$ and $\sP(\Om_2)$, respectively:
\begin{align*}
  Q_N(\omega) :=\oo{N} \sum_{i=1}^N \delta_{Y^{i,N}(\omega)},
  \qquad
  R_N(\omega) := \oo{N(N-1)} \sum_{i \neq j} \delta_{\bar
  Y^{i,j,N}(\omega)}.
\end{align*}
Conditionally on the system, $R_N$ is the law of
$(Y^{I,N}, Y^{J,N}, \ell^{I,J,N})$, where the ordered pair $(I,J)$ of
distinct particles is chosen uniformly at random, independently of
the system. We call a realization $R_N(\omega)$ of $R_N$ the \define{quenched} law of a random pair.
Likewise, $Q_N(\omega)$ is the quenched law of a random particle. 
With
$\pi_1(y, \bar y, l) := y$ and
$\pi_{1,2}(y, \bar y, l) := (y, \bar y)$, we have
$Q_N = \pushf{\pi_1}{R_N}$. If we, additionally, average 
over the system, we get the
\define{annealed} law, i.e., the mean measure $\ee{R_N}$.
By the exchangeability established in Section~\ref{ss-bset}, it is the law of
$\bar Y^{1,2,N}$:
\begin{equation}\label{mean}
  \ee{\inner{R_N}{\vp}} = \ee{\vp(\bar Y^{1,2,N})}
  \text{ for bounded measurable } \vp \colon \Om_2 \to \R .
\end{equation}

The jump measure of $x \in D([0,\infty),\R)$ is
\begin{align*}
  \Gamma(x) := \sum_{s > 0 \colon \Delta x_s \neq 0}
  \delta_{(s, \Delta x_s)} \in \sM_p((0,\infty) \times E) .
\end{align*}
It is locally finite, since $x$ has finitely many jumps of size at
least $\eps$ on $[0,T]$, for all $\eps, T > 0$. Since $\mu^{i,N}$ has
unit atoms, at most one per time, and $B^{i,N}$ is continuous,
\eqref{X} gives
\begin{equation}\label{jumpmeas}
  \mu^{i,N} = \Gamma(X^{i,N}) \quad \text{a.s.}
\end{equation}
We record here a property of the jump functional $\Gamma$ 
for later use in Section \ref{s-ident} below. 
\begin{lemma}\label{l-jump}
  \begin{enumerate}[(a)]
  \item\label{i-jump-cont} The map $\Gamma \colon D([0,\infty),\R) \to \sM_p((0,\infty) \times
  E)$ is continuous.
  \item\label{i-jump-sum} For $h \in C_c(E)$ and $t > 0$, the map
    $x \mapsto \sum_{0 < s \leq t} h(\Delta x_s)$ is continuous at
    every $x \in D([0,\infty),\R)$ with $\Delta x_t = 0$.
  \end{enumerate}
\end{lemma}

\begin{proof}
  \ref{i-jump-cont} Let $x_n \to x$. There are increasing homeomorphisms $\rho_n$ of
  $[0,\infty)$ with $\sup_s \abs*{\rho_n(s) - s} \to 0$ and
  $\tilde x_n := x_n \circ \rho_n \to x$ locally uniformly
  \cite[Theorem VI.1.14]{JacShi03}. Let $f \in C_c((0,\infty) \times E)$
  vanish unless $s \leq T$ and $\abs*{z} \geq \eps$. Then
  $\int f \, \diffd\Gamma(x_n) = \sum_s f(\rho_n(s), \Delta \tilde x_n(s))$,
  the sum over the jump times of $\tilde x_n$. For large $n$,
  $\rho_n(s) \leq T$ implies $s \leq T + 1$, and
  $\abs*{\Delta \tilde x_n(s) - \Delta x(s)} \leq 2 \sup_{r \leq T+1}
  \abs*{\tilde x_n(r) - x(r)} < \eps/2$ for $s \leq T + 1$. So only the
  finitely many $s \leq T + 1$ with $\abs*{\Delta x(s)} \geq \eps/2$
  contribute, and at each of them $\rho_n(s) \to s$ and
  $\Delta \tilde x_n(s) \to \Delta x(s) \neq 0$. Since $f$ is
  continuous, $\int f \, \diffd\Gamma(x_n) \to \int f \, \diffd\Gamma(x)$.

  \ref{i-jump-sum} Since $h$ vanishes near $0$, the map
  $x \mapsto x^h$, $x^h(r) := \sum_{0 < s \leq r} h(\Delta x_s)$, is
  continuous from $D([0,\infty),\R)$ into itself
  \cite[Corollary VI.2.8]{JacShi03}. Every jump of $x^h$ is a jump
  of $x$, so $\Delta x^h_t = 0$. Evaluation at $t$ is continuous at
  such paths \cite[VI.2.3]{JacShi03}.
\end{proof}

\subsection{Tightness}\label{ss-tightness}

We prove the necessary tightness results next. We use two
standard tools, Aldous's criterion and a tightness criterion for
random measures, and explain below how they work in our setting.

  At the level of one pair of
particles, we show that the laws of $\bar Y^{1,2,N}$ are tight in
$\sP(\Om_2)$. By exchangeability, it suffices to treat the components
$X^{1,N}$, $B^{1,N}$, $\ld^{1,N}$ and $\ell^{1,2,N}$. For each of
them we use Aldous's criterion \cite[Theorem VI.4.5]{JacShi03}. It
asks for compact containment on finite horizons, and for increments
over stopping intervals of length at most $\th$ that are small in
probability, uniformly in $N$, as $\th \to 0$. Compact containment
comes from the moment bounds of Section~\ref{s-bounds}, and the
increments are controlled by Proposition~\ref{p-aldous}.

At the level of the empirical measures, we use a standard 
criterion for random measures as presented in 
\cite[Theorem 4.10]{Kal17}. It says
that a family $\Xi$ of almost surely finite random measures on a
Polish space $S$ is tight for the weak topology if and only if
\begin{enumerate}[(i)]
\item\label{i-kal-mass} the total masses $\xi(S)$, $\xi \in \Xi$,
  are tight in $[0,\infty)$, and
\item\label{i-kal-tail} $\inf_K \sup_{\xi \in \Xi}
  \ee*{\xi(S \setminus K) \wedge 1} = 0$, where the infimum is over
  compact $K \subseteq S$.
\end{enumerate}
For random probability measures, condition~\ref{i-kal-mass} is
automatic, since $\xi(S) = 1$. Condition~\ref{i-kal-tail} becomes
$\inf_K \sup_{\xi} \ee*{\xi(S \setminus K)} = 0$, which is tightness
of the mean measures $A \mapsto \ee*{\xi(A)}$. By \eqref{mean}, the
mean measure of $R_N$ is the law of $\bar Y^{1,2,N}$. So the compact sets
that work for the laws of $\bar Y^{1,2,N}$ work for all $R_N$ at once, and
tightness of the laws of $R_N$ reduces to the first level. In the
language of Section~\ref{ss-paths}, tightness of the quenched laws
of a random pair reduces to tightness of the annealed law. 
\begin{proposition}\label{p-tight}
  \begin{enumerate}[(a)]
  \item\label{i-tight-U} The sequence $\seq[N]{\bar Y^{1,2,N}}$ is
    tight in $\Om_2$. In particular, $\seq[N]{Y^{1,N}}$ is tight in
    $\Om_1$.
  \item\label{i-tight-R} The sequences $\seq[N]{R_N}$ and
    $\seq[N]{Q_N}$ are tight in $\sP(\Om_2)$ and $\sP(\Om_1)$.
  \item\label{i-tight-prod} If $\mathfrak R \in \sP(\sP(\Om_2))$ is
    the weak limit of the laws of a subsequence of $\seq[N]{R_N}$, then, for
    $\mathfrak R$-almost every $\rho$,
    \begin{align*}
      \pushf{(\pi_{1,2})}{\rho}
      = \pushf{\pi_1}{\rho} \otimes \pushf{\pi_1}{\rho} .
    \end{align*}
    That is, the first two coordinates are independent under
    $\rho$, each with law $\pushf{\pi_1}{\rho}$.
  \end{enumerate}
\end{proposition}

\begin{proof}
  \ref{i-tight-U} A family of laws on a product of Polish spaces is
  tight if and only if the families of its marginals are (take
  products of $\eps$-supporting compact sets for one direction,
  and continuous projections for the other). By
  exchangeability, $Y^{2,N}$ has the law of $Y^{1,N}$, so it
  suffices to treat $X^{1,N}$, $B^{1,N}$, $\ld^{1,N}$ and
  $\ell^{1,2,N}$. The law of $B^{1,N}$ does not depend on $N$.

  For $X^{1,N}$ we use Aldous's criterion
  \cite[Theorem VI.4.5]{JacShi03}.
  By \eqref{X}, for $T \geq 0$,
  \begin{align*}
    \sup_{t \leq T} \abs*{X^{1,N}_t}
    \leq \abs*{X^{1,N}_0} + \sup_{t \leq T} \abs*{B^{1,N}_t}
    + \sum_{s \leq T} \abs*{\Delta X^{1,N}_s} .
  \end{align*}
  The laws of $X^{1,N}_0$ converge by Assumption~\ref{a-std}\ref{i-chaos},
  so they are tight. For $m \in \N$ and $b > 0$, let
  \begin{align*}
    A_{m,b} := \set*{J^{1,N}_T \leq m} \cap
    \set*{\mu^{1,N}\prn*{(0,T] \times \sets{z}{\abs{z} > b}} = 0}
  \end{align*}
  be the event that $X^{1,N}$ has at most $m$ jumps in $(0,T]$, all
  of size at most $b$. On $A_{m,b}$, the last sum in the display is
  at most $mb$. By Markov's inequality and Lemma~\ref{l-jumpmom}\ref{i-bdnu},
  \begin{align*}
    \pp*{A_{m,b}^c}
    \leq \frac{\ee*{J^{1,N}_T}}{m}
    + \ee*{\mu^{1,N}\prn*{(0,T] \times \sets{z}{\abs{z} > b}}}
    \leq \frac{C_T}{m} + C_T \, \nu\prn*{\sets{z}{\abs{z} > b}} .
  \end{align*}
  The right side does not depend on $N$. Given $\eps > 0$, choose
  $m$ with $C_T/m < \eps$, and then $b$ with
  $C_T \, \nu(\sets{z}{\abs{z} > b}) < \eps$, which is possible since
  $\nu$ is a probability measure on $E$. Choose also $r$ with
  $\pp*{\abs*{X^{1,N}_0} > r} < \eps$ for all $N \geq 2$ and
  $\pp*{\sup_{t \leq T} \abs*{B^{1,N}_t} > r} < \eps$. The latter
  probability does not depend on $N$. Then
  \begin{align*}
    \pp*{\sup_{t \leq T} \abs*{X^{1,N}_t} > 2r + mb} < 4 \eps
    \quad \text{for all } N \geq 2 .
  \end{align*}
  Since $\eps$ is arbitrary, the sequence
  $\seq[N]{\sup_{t \leq T} \abs*{X^{1,N}_t}}$ is tight. This is
  the compact containment condition in Aldous's criterion. For the
  increments, let
  $\sg \leq \tau \leq \sg + \th$ be stopping times with $\sg \leq T$ and
  $\th \in [0,1]$, and let $\eps \in (0,1]$. Since
  $\abs*{X^{1,N}_{\tau} - X^{1,N}_{\sg}} \geq \eps$ implies
  $1 \wedge \abs*{X^{1,N}_{\tau} - X^{1,N}_{\sg}} \geq \eps$, Markov's
  inequality and Proposition~\ref{p-aldous} give
  \begin{align*}
    \pp*{\abs*{X^{1,N}_{\tau} - X^{1,N}_{\sg}} \geq \eps}
    \leq \frac{C_T \sqrt{\th}}{\eps},
  \end{align*}
  uniformly in $N$, $\sg$ and $\tau$. 
  This tends to $0$ as $\th \to 0$.

  The same argument applies to $\ld^{1,N}$ and $\ell^{1,2,N}$, with
  compact containment from Corollary~\ref{c-mom}\ref{i-mom-sup} for
  $\ld^{1,N}$. For $\ell^{1,2,N}$ it follows from
  $\ee*{\ell^{1,2,N}_T} \leq C_T$, which Proposition~\ref{p-aldous}
  gives after splitting $[0,T]$ into intervals of length at most $1$
  (note that $\sup_{t \leq T} \ell^{1,2}_t =
  \ell^{1,2}_T$). Their paths are continuous. The space of continuous
  paths is closed in the Skorokhod space, and the Skorokhod topology
  induces the locally uniform topology on it
  \cite[Theorem VI.1.14, Proposition VI.1.17]{JacShi03}. So tightness in the Skorokhod space
  gives tightness in $C([0,\infty),\R)$, and in its closed subsets
  $C([0,\infty),[0,\infty))$ and $C_\uparrow([0,\infty))$.

  \ref{i-tight-R} By \eqref{mean} and \ref{i-tight-U}, the mean
  measures of the $R_N$ are tight. By the discussion before the
  proposition, the laws of the $R_N$ are tight. The map
  $\rho \mapsto \pushf{\pi_1}{\rho}$ is continuous, and
  $Q_N = \pushf{\pi_1}{R_N}$.

  \ref{i-tight-prod} Splitting $Q_N \otimes Q_N$ into its diagonal
  and off-diagonal parts gives
  \begin{align*}
    \pushf{(\pi_{1,2})}{R_N} - Q_N \otimes Q_N
    = \oo{N-1} \, Q_N \otimes Q_N
    - \oo{N(N-1)} \sum_{i=1}^N \delta_{(Y^{i,N}, Y^{i,N})} .
  \end{align*}
  Both terms on the right have total mass $1/(N-1)$. Let $d$ be the
  bounded Lipschitz metric on $\sP(\Om_1 \times \Om_1)$. That is, we
  fix a metric $d_0$ on $\Om_1 \times \Om_1$ that induces its topology,
  and put
  \begin{align*}
    d(m, m') := \sup \sets*{\abs*{\tint f \, \diffd m - \tint f \, \diffd m'}}
    {\abs{f} \leq 1,\ \abs{f(u) - f(v)} \leq d_0(u,v)
    \text{ for all } u, v},
  \end{align*}
  for $m, m' \in \sP(\Om_1 \times \Om_1)$, where $f$ ranges over
  functions on $\Om_1 \times \Om_1$. It metrizes weak
  convergence and is dominated by the total variation norm. The maps
  $\rho \mapsto \pushf{(\pi_{1,2})}{\rho}$ and
  $\rho \mapsto \pushf{\pi_1}{\rho} \otimes \pushf{\pi_1}{\rho}$ are
  continuous, so
  $h(\rho) := d(\pushf{(\pi_{1,2})}{\rho}, \pushf{\pi_1}{\rho}
  \otimes \pushf{\pi_1}{\rho})$ is continuous on $\sP(\Om_2)$, and
  $h(R_N) \leq 2/(N-1)$. Since $h$ is also bounded,
  $\int h \, \diffd\mathfrak R = \lim \ee*{h(R_N)} = 0$ along the
  subsequence, so $h = 0$, $\mathfrak R$-almost everywhere.
\end{proof}

\section{The nonlinear system}\label{s-nonlin}

By \eqref{ld}, the intensity of particle $i$
is driven by $(N-1) \gm_N \to \kp$ multiplied by the average
$\oo{N-1} \sum_{j \neq i} \ell^{i,j,N}$ of its collision local
times. As $N \to \infty$, we expect the particles to
decouple, and this average to approach the conditional
expectation, given particle $i$, of its collision local time with an
independent particle of the same law. By Lemma~\ref{l-condlt} below,
this conditional expectation is $\int_0^t u(s, X^i_s) \, \ds$, where
$u(s,\cdot)$ is the density of the position at time $s$. This leads to a
one-particle system whose intensity is driven by the density of its
own position. Its dynamics depend, through $u$, on the law of
the solution itself, and not only on its current state. So the
evolution equation for this law is nonlinear. As is customary for
such McKean--Vlasov type systems, we call it the \define{nonlinear
system}. We define it and prove uniqueness. Existence will
follow from the convergence of the particle system.

\subsection{Definition}\label{ss-nldef}
Let $(X, B, \ld)$ be the coordinate process on the path space
$\Om_1$ of Section~\ref{ss-paths}, and let $\mu := \Gamma(X)$ be the
jump measure of $X$. Let $\bF^0 = (\sF^0_t)_{t \geq 0}$,
where $\sF^0_t$ is generated by $X_r$, $B_r$ and $\ld_r$ for
$r \leq t$. For
$P \in \sP(\Om_1)$, let $\bF^P$ be the usual $P$-augmentation of
$\bF^0$, and write $\ee^{P}*{\cdot}$ for the expectation under $P$.

\begin{definition}\label{d-nl}
  A \define{solution of the nonlinear system} is a law
  $P \in \sP(\Om_1)$ with the following properties.
  \begin{enumerate}[(a)]
  \item\label{i-nl-init} Under $P$, $(X_0, \ld_0)$ has law $\eta_0$.
  \item\label{i-nl-bm} $B$ is an $\bF^P$-Brownian motion.
  \item\label{i-nl-graph} Under $P$, $X$ is piecewise Brownian
    with variance $1$ in $\bF^P$, with $W = B$.
  \item\label{i-nl-comp} The $\bF^P$-compensator of $\mu$ is
    $\ld_s \, \ds \, \nu(\dz)$.
  \item\label{i-nl-mom} $\sup_{s \leq T} \ee^{P}*{\ld_s} < \infty$ for
    every $T \geq 0$.
  \item\label{i-nl-feed} For every $t > 0$, the law of $X_t$
    under $P$ has a density $u(t,\cdot)$, with $u(\cdot,\cdot)$
    jointly measurable, such that 
    $P$-almost surely,
    for all $t \geq 0$,
    \begin{equation}\label{feed}
      \ld_t = e^{-t} \ld_0
      + \kp \int_0^t e^{-(t-s)} u(s, X_s) \, \ds.
    \end{equation}
  \end{enumerate}
\end{definition}
\subsection{Densities}\label{ss-dens}
Let $p_r(x) := (2\pi r)^{-1/2} e^{-x^2/(2r)}$ be the heat kernel, and
put
\begin{align*}
  q_r(x) := (p_r*\nu) (x) - p_r(x), \qquad r > 0,
\end{align*}
where, for a function $f$ and a finite measure $m$ on $\R$, 
$(f * m)(x) := \int f(x - y) \, m(\dy)$ denotes the convolution of
$f$ and $m$. Let $\vth_0$ 
be 
the law of the initial position, i.e., 
the image of
$\eta_0$ under $(x,y) \mapsto x$. 
Since $\nu$ is a probability measure, both $p_r$ and $p_r * \nu$
take values in $[0,(2\pi r)^{-1/2}]$. Thus
$\norm{q_r}_\infty \leq \norm{p_r}_\infty = (2\pi r)^{-1/2}$,
which gives the elementary bound
\begin{equation}\label{heat}
  \norm{p_r}_\infty + \norm{q_r}_\infty \leq r^{-1/2} .
\end{equation}

\begin{lemma}\label{l-dens}
  Let $P \in \sP(\Om_1)$ satisfy
  \ref{i-nl-init}--\ref{i-nl-mom} of Definition~\ref{d-nl}. For
  $t \geq 0$, let $w_t$ be the finite measure
  $w_t(A) := \ee^{P}*{\ld_t \ind{A}(X_t)}$. Then the function
  \begin{equation}\label{mild}
    u(t,\cdot) := p_t * \vth_0 + \int_0^t q_{t-s} * w_s \, \ds
  \end{equation}
  is continuous on $(0,\infty) \times \R$, and $u(t,\cdot)$ is a density of
  the law of $X_t$ for every $t > 0$.
\end{lemma}

\begin{proof}
  Fix $t > 0$ and $\psi \in C_c^\infty(\R)$. The
  function $\vp(s,x) := p_{t-s} * \psi(x)$ solves
  $\partial_s \vp + \ot \partial_x^2 \vp = 0$ on $[0,t) \times \R$.
  Apply \itos\ formula to $\vp(s, X_s)$ on $[0, t - \dl]$, take
  expectations, and let $\dl \downarrow 0$. By \ref{i-nl-graph},
  \ref{i-nl-comp} and \ref{i-nl-mom}, the jumps contribute their
  compensator, and
  \begin{multline*}
    \ee^{P}{\psi(X_t)} =  \int (p_t * \psi) \, \diffd\vth_0 +
    \\ +\int_0^t \int \prn{\int_E (p_{t-s} * \psi)(x + z) \, \nu(\dz)
    - (p_{t-s} * \psi)(x)} w_s(\dx) \, \ds .
  \end{multline*}
  By Fubini's theorem and since $p_r$ is even,
  $\int (p_r * \psi) \, \diffd m = \int \psi \, (p_r * m) \, \dx$ for every
  finite measure $m$, and likewise for the shifted term. Hence the
  right side equals $\int \psi(x) \, u(t,x) \, \dx$.
  By \eqref{heat}, $\norm*{q_r * w_s}_\infty \leq r^{-1/2} \,
  \ee^{P}*{\ld_s}$, and \ref{i-nl-mom} and dominated convergence give
  continuity.
\end{proof}
\begin{remark}
\label{r-mild}
\begin{enumerate}[(1)]
  \item\label{i-mild-fwd} Formula \eqref{mild} is the mild form, by Duhamel's principle, of
the forward equation
\begin{equation}\label{fwd}
  \partial_t u = \ot \partial_x^2 u + \nu * w_t - w_t,
  \qquad u(0,\cdot) = \vth_0,
\end{equation}
in the sense of distributions, where
$(\nu * w_t)(A) := \int_E w_t(A - z) \, \nu(\dz)$.
\item\label{i-mild-version} If $P$ is a solution of the nonlinear system, we always take $u$ in
\eqref{feed} to be the continuous version \eqref{mild} of
Lemma~\ref{l-dens}, here and in the rest of the paper. We can do
this without loss of generality: two jointly measurable versions
of the density agree, for every $s$, almost everywhere in $x$, and
$X_s$ has a density, so they give the same integral
$\int_0^t u(s, X_s) \, \ds$, $P$-almost surely, simultaneously for all
$t$.
\end{enumerate}
\end{remark}
\begin{lemma}\label{l-densbd}
  Let $P$ be a solution of the nonlinear system. Then, for every
  $T \geq 0$,
  \begin{align*}
    \sup_{0 < t \leq T} \sqrt{t} \, \norm{u(t,\cdot)}_\infty \leq C_T
    \quad\text{and}\quad
    \sup_{t \leq T} \ld_t \leq \ld_0 + C_T \quad P\text{-a.s.}
  \end{align*}
\end{lemma}

\begin{proof}
  Fix $T$, and let
  $g(t) := \sup_{0 < s \leq t} \sqrt{s} \, \norm*{u(s,\cdot)}_\infty$, which is
  finite by Lemma~\ref{l-dens}, \eqref{heat} and \ref{i-nl-mom}. By
  \eqref{feed},
  $\ld_s \leq \ld_0 + \kp \int_0^s r^{-1/2} g(r) \, \dr$, so
  $w_s(\R) = \ee^{P}*{\ld_s}$ is at most the same expression with
  $\ee^{P}*{\ld_0}$ in place of $\ld_0$. Here
  $\ee^{P}*{\ld_0} = \int y \, \diffd\eta_0 \leq x_G + 1 \leq C_T$, by
  \eqref{G0}. Insert
  this into \eqref{mild}, use \eqref{heat}, and exchange the order of
  integration. For $0 < t \leq T$,
  \begin{align*}
    \sqrt{t} \, \norm{u(t,\cdot)}_\infty \leq 1 + 2 T \, \ee^{P}{\ld_0}
    + 2 \kp T \int_0^t r^{-1/2} g(r) \, \dr .
  \end{align*}
  The right side is nondecreasing in $t$, so $g$ satisfies the same
  bound. Gronwall's lemma \cite[Appendix~5,
  Theorem~5.1]{EthKur86}, with the integrable kernel $r^{-1/2}$,
  gives $g(T) \leq C_T$. Then $\sup_{t \leq T} \ld_t \leq
  \ld_0 + 2 \kp \sqrt{T} \, C_T$.
\end{proof}

\subsection{Conditional collision local time}\label{ss-condlt}
Following Sznitman's approach to collision-local-time interactions
\cite{Szn86}, we identify the averaged collision input through
the conditional local time of two independent copies.
The next lemma identifies the driving term in \eqref{feed} as a
conditional collision local time. By Lemma~\ref{l-densbd}, its
assumptions hold for every solution of the nonlinear system. We do
not assume \ref{i-nl-feed}, since the lemma will be used in
Section~\ref{s-ident} below before \ref{i-nl-feed} is known.

\begin{lemma}\label{l-condlt}
  Let
  \begin{enumerate}[(i)]
  \item $P$ satisfy \ref{i-nl-init}--\ref{i-nl-mom} of
  Definition~\ref{d-nl},
  and $u$ be the continuous density of Lemma~\ref{l-dens},
  \item $Y^k = (X^k, B^k, \ld^k)$, $k = 1, 2$, be independent with law
  $P$, and $\mu^k := \Gamma(X^k)$,
  \item $\ell := \ell(0)$ be
  defined by \eqref{tanaka}, with $c = 2$,
  $Z := X^1 - X^2$, $W := B^1 - B^2$
  and the usual augmentation
  of the filtration generated by $Y^1$ and $Y^2$.
  \end{enumerate}
  Then, for every
  $t \geq 0$,
  \begin{equation}\label{condlt}
    \ee{\ell_t \giv Y^1} = \int_0^t u(s, X^1_s) \, \ds
    \quad \text{a.s.}
  \end{equation}
\end{lemma}

\begin{proof}
  For $\eps \in (0,1]$, put
  \begin{align*}
    f_\eps := \frac{1}{2\eps} \int_0^t \inds{\abs*{Z_s} < \eps} \, \ds .
  \end{align*}
  The idea is to show that $f_\eps \to \ell_t$ in $\bL^1$, and that
  $\ee*{f_\eps \giv Y^1} \to \int_0^t u(s, X^1_s) \, \ds$ in $\bL^1$,
  as $\eps \downarrow 0$. This implies that 
  $\ee*{f_\eps \giv Y^1} \to \ee*{\ell_t \giv Y^1}$ in $\bL^1$, and
  \eqref{condlt} follows.

  \emph{Step 1: $f_\eps \to \ell_t$ in $\bL^1$.} In the filtration of
  the lemma, $B^1$ and $B^2$ are independent Brownian motions, and $Z$
  has finitely many jumps on bounded intervals. So $Z$ is piecewise
  Brownian with variance $2$, with $V := Z - Z_0 - W$. Since
  the local time $\ell_t(a)$ of $Z$ at $a$ is an occupation density,
  \begin{align*}
    f_\eps = \frac{1}{2\eps} \int_{-\eps}^{\eps} \ell_t(a) \, \da
    = \int \prn{\ell_t(a) - \ell_0(a)} \rho_\eps(\da),
  \end{align*}
  where $\rho_\eps$ is the uniform distribution on $(-\eps,\eps)$. As
  $a \mapsto \ell_t(a)$ is continuous, $f_\eps \to \ell_t$ almost
  surely. For uniform integrability, we use
  Lemma~\ref{l-pieces}\ref{i-pieces-ui} with $\sg = 0$,
  $\tau = \th = t$, $c = 2$ and $\rho = \rho_\eps$. The random
  variable $K$ of \eqref{num-pieces} is at most
  $1 + \mu^1((0,t] \times E) + \mu^2((0,t] \times E)$, which is
  integrable by \ref{i-nl-comp} and \ref{i-nl-mom}, so that
  $\vp(M) := \ee*{(K - M)^+} \to 0$ as $M \to \infty$ by dominated
  convergence. The same $Z$, $\sg$, $\tau$ and $K$ serve for all $\eps$.
  Thus, $(f_\eps)_{\eps \in (0,1]}$
  is uniformly integrable, and $f_\eps \to \ell_t$ in $\bL^1$.

  \emph{Step 2: the conditional expectations.} Write
  $f_\eps = \int_0^t h_\eps(X^1_s - X^2_s) \, \ds$ with
  $h_\eps := (2\eps)^{-1} \ind{(-\eps,\eps)}$, and put
  \begin{align*}
    u^\eps(s, x) := \int h_\eps(x - y) \, u(s, y) \, \dy
    = \frac{1}{2\eps} \int_{x - \eps}^{x + \eps} u(s, y) \, \dy .
  \end{align*}
  For each $s$, $X^1_s$ is a function of $Y^1$ and $X^2_s$ is
  independent of $Y^1$ with density $u(s,\cdot)$, so Fubini's theorem
  gives $\ee*{f_\eps \giv Y^1} = \int_0^t u^\eps(s, X^1_s) \, \ds$
  almost surely.
  Since $u$ is continuous by Lemma~\ref{l-dens}, $u^\eps \to u$
  pointwise on $(0,t] \times \R$. Moreover,
  by \eqref{mild},
  \eqref{heat} and \ref{i-nl-mom},  we have
\begin{align*}
  u^\eps(s, x) \leq \norm{u(s, \cdot)}_\infty \leq
    s^{-1/2}\prn{1+ 2t
  \textstyle\sup_{r \leq t} \ee^{P}*{\ld_r}}
 \efor 0 < s \leq t.
\end{align*}
  The function $s^{-1/2}$ is integrable on $(0,t]$ so, by
  dominated convergence,
  $\int_0^t u^\eps(s, X^1_s) \, \ds \to \int_0^t u(s, X^1_s) \, \ds$
  almost surely and, since the bound is deterministic, in $\bL^1$.
\end{proof}

\begin{remark}\label{r-meanfield}
  Let $P$ be a solution of the nonlinear system. Put $A_t := \ee*{\ell_t \giv Y^1}$, the conditional collision local
  time of the two copies. By \eqref{condlt}, \eqref{feed} for the
  first copy reads
  \begin{align*}
    \ld^1_t = e^{-t} \ld^1_0 + \kp \int_0^t e^{-(t-s)} \, \diffd A_s .
  \end{align*}
  This is \eqref{ld} with $\gm_N \sum_{j \neq i} \ell^{i,j}$ replaced
  by $\kp A$, the natural candidate for its mean-field limit. Indeed,
  $(N-1) \gm_N \to \kp$, and the average of $\ell^{i,j}$ over $j$
  should approach the conditional collision local time with an
  independent particle. This is made rigorous in
  Section~\ref{s-ident} below. 
\end{remark}

\subsection{Uniqueness}\label{ss-nlunique}
Let $(\Om, \sF, \bF, \bP)$ support an $\sF_0$-measurable random
vector $(X_0, \ld_0)$ with law
$\eta_0$, an $\bF$-Brownian motion $B$, and an $\bF$-Poisson random
measure $\Pi$ on $(0,\infty) \times E \times (0,\infty)$ with
intensity $\ds \, \nu(\dz) \, \dy$, all independent of each other. 
We call the vector $( X_0,\ld_0, B, \Pi)$ a \define{noise}. 

For a Borel function
$v \colon (0,\infty) \times \R \to [0,\infty)$ with
$\int_0^T \norm*{v(s,\cdot)}_\infty \, \ds < \infty$ for each $T$,
and noise $(X_0, \ld_0, B, \Pi)$ consider the
system
\begin{equation}\label{frozen}
  \begin{aligned}
    X_t &= X_0 + B_t + \tint_{(0,t] \times E} z \, \mu(\ds,\dz), \\
    \mu(\ds,\dz) & = \tint_{(0,\infty)} \inds{y \leq \ld_s}
    \, \Pi(\ds,\dz,\dy), \\
    \ld_t &= e^{-t} \ld_0 + \kp \tint_0^t e^{-(t-s)} v(s, X_s) \, \ds .
  \end{aligned}
\end{equation}
In any solution of \eqref{frozen}, $\ld$ is bounded on bounded
intervals, since $\int_0^T \norm*{v(s,\cdot)}_\infty \, \ds < \infty$.
The measure $\mu$ retains only the atoms of $\Pi$ that lie below the
graph of $\ld$, and almost surely $\Pi$ has finitely many atoms below
any level on bounded time intervals. Therefore $\mu$ has only finitely many
atoms on $(0,T] \times E$ for each $T \ge 0$, a.s.
Therefore, the construction in the
proof of Proposition~\ref{p-exist}, with one particle and with
\eqref{ld} replaced by the last equation in \eqref{frozen}, gives a
solution $(X, \ld, \mu)$ which is, moreover, of the form
$(X,\ld,\mu) = \Phi_v(X_0, \ld_0, B, \Pi)$ for a measurable map
$\Phi_v$.
Solutions of \eqref{frozen} are unique in the following sense.
Every solution driven by a noise, on any filtered probability space,
equals $\Phi_v$ applied to that noise. To see this, argue by
induction over the atoms of $\Pi$ that contribute to $\mu$, taken in
the order of their time coordinates, as in
Remark~\ref{r-thin-unique}.

The last equation in \eqref{frozen} is \eqref{feed} with $u$
replaced by $v$, so that the nonlinear
system can be rephrased, essentially, as a fixed point problem: $u$ is
the density of the $X$-component of  $\Phi_u$. 
We use this to prove uniqueness below by running 
$\Phi_u$ and $\Phi_{\bar u}$ with the same noise and closing the
loop between the distance of the outputs and that of $u$ and
$\bar u$, through \eqref{mild}, by Gronwall's lemma.

\begin{lemma}\label{l-thin}
  Let $P$ be a solution of the nonlinear system, with density $u$.
  Given a noise $(X_0, \ld_0, B,\Pi)$ on
  $(\Om, \sF, \bF, \bP)$ as above, define
  \begin{align*}
  (\hX, \hat\ld, \hat\mu) := \Phi_u(X_0, \ld_0, B, \Pi) .
  \end{align*}
  Then $(\hX, B, \hat\ld)$ has law $P$.
\end{lemma}

The proof is given in Appendix~\ref{ap-thin}.

\begin{proposition}\label{p-nlunique}
  There is at most one solution of the nonlinear system.
\end{proposition}

\begin{proof}
  Let $P$ and $\bar P$ be solutions, with densities $u$ and $\bar u$.
  Take a noise $(X_0, \ld_0, B, \Pi)$ as above, and put
  $(X, \ld, \mu) := \Phi_u(X_0, \ld_0, B, \Pi)$ as well as 
  $(\bar X, \bar\ld, \bar\mu) := \Phi_{\bar u}(X_0, \ld_0, B,
  \Pi)$, coupling the two solutions together. Indeed,  
  by Lemma~\ref{l-thin}, $(X, B, \ld)$ has law $P$ and
  $(\bar X, B, \bar\ld)$ has law $\bar P$. 

  For a fixed $T$
  Lemma~\ref{l-densbd} implies that
  $\norm*{u(s,\cdot)}_\infty + \norm*{\bar u(s,\cdot)}_\infty \leq C_T s^{-1/2}$ for
  $0 < s \leq T$, and $\ld, \bar\ld \leq \ld_0 + C_T$ on $[0,T]$. 
  For
  $t \leq T$, let
  \begin{align*}
    D_t := \ee{(1 + \ld_0) \inds{X_s \neq \bar X_s
    \text{ for some } s \leq t}}
    + \sup_{s \leq t} \ee{(1 + \ld_0) \abs*{\ld_s - \bar\ld_s}} .
  \end{align*}
  It is finite, since $\ee{\ld_0} < \infty$ by \eqref{G0}, and
  since $\ld$ and $\bar\ld$ share the same $\ld_0$, so that, by
  \eqref{feed},
  $\abs*{\ld_s - \bar\ld_s} \leq \kp \int_0^s \prn{\norm*{u(r,\cdot)}_\infty
  + \norm*{\bar u(r,\cdot)}_\infty} \dr \leq C_T$ for $s \leq T$.

  \emph{Positions.} The positions agree up to the first atom of $\Pi$
  ``accepted'' by exactly one of $\ld$ and $\bar\ld$. Both are continuous
  and adapted, and $\nu(E) = 1$, so the number of such atoms in
  $(0,t]$ has compensator $\int_0^t \abs*{\ld_s - \bar\ld_s} \, \ds$.
  Since $1 + \ld_0$ is $\sF_0$-measurable,
  \begin{align*}
    \ee{(1 + \ld_0) \inds{X_s \neq \bar X_s \text{ for some }
    s \leq t}}
    \leq \int_0^t \ee{(1 + \ld_0) \abs*{\ld_s - \bar\ld_s}} \, \ds
    \leq \int_0^t D_s \, \ds .
  \end{align*}

  \emph{Densities.} We use the representation and notation of
  subsection \ref{ss-dens} and Lemma~\ref{l-dens}. In particular 
  let $w$ and $\bar w$ be as defined there, using $X$ and $\bar X$,
  respectively, so that 
\begin{align*}
  (q_r * w_s)(x) = \ee*{\ld_s \, q_r(x - X_s)} \eand
  (q_r * \bar w_s)(x) = \ee*{\bar\ld_s \, q_r(x - \bar X_s)}.
\end{align*}
  By \eqref{heat},
  $\abs*{q_r(a) - q_r(b)} \leq 2 r^{-1/2} \inds*{a \neq b}$. Since
  $\bar\ld_s \leq (1 + C_T)(1 + \ld_0)$, we have
  \begin{align*}
    \norm{q_r * w_s - q_r * \bar w_s}_\infty
    &{}\leq r^{-1/2} \, \ee{\abs*{\ld_s - \bar\ld_s}}
    + 2 r^{-1/2} \, \ee{\bar\ld_s \inds{X_s \neq \bar X_s}}
 \\ &\leq C_T r^{-1/2} D_s .
  \end{align*}
  The initial terms $p_t * \vth_0$ of \eqref{mild} do not depend on $w$ or $\bar
  w$, so they cancel in the difference. Hence
  $\norm*{u(t,\cdot) - \bar u(t,\cdot)}_\infty \leq C_T \int_0^t (t-s)^{-1/2} D_s \,
  \ds$, and by Fubini's theorem
\begin{align*}
  \int_0^t \norm*{u(s,\cdot) - \bar u(s,\cdot)}_\infty \, \ds
  &\leq  C_T \int_0^t 2 (t-r)^{1/2} D_r \, \dr
  \leq C_T \int_0^t D_r \, \dr .
\end{align*}

  \emph{Intensities.} Since
  $\abs*{u(s, X_s) - \bar u(s, \bar X_s)} \leq 2 \norm*{u(s,\cdot)}_\infty
  \inds*{X_s \neq \bar X_s} + \norm*{u(s,\cdot) - \bar u(s,\cdot)}_\infty$,
  equation \eqref{feed} and $\ee*{1 + \ld_0} \leq C_T$ give
  \begin{align*}
    \ee{(1 + \ld_0) \abs*{\ld_t - \bar\ld_t}}
    &\leq C_T \int_0^t \prn{s^{-1/2} D_s
    + \norm{u(s,\cdot) - \bar u(s,\cdot)}_\infty} \ds
 \\ & \leq C_T \int_0^t \prn{1 + s^{-1/2}} D_s \, \ds .
  \end{align*}

  Combining the three bounds,
  $D_t \leq C_T \int_0^t (1 + s^{-1/2}) D_s \, \ds$ for $t \leq T$.
  The kernel is integrable, so Gronwall's lemma
  gives $D = 0$ on
  $[0,T]$. Hence $X = \bar X$ and $\ld = \bar\ld$ on $[0,T]$ and,
  as $T$ is arbitrary, $P = \bar P$.
\end{proof}

\section{Identification of the limit and propagation of chaos}\label{s-ident}\label{s-chaos}
We identify the limit points of the laws of $R_N$. Each of them
is concentrated on laws of a pair of independent processes, whose
common law solves the nonlinear system, together with their collision
local time (Theorem~\ref{t-ident}). Uniqueness of the nonlinear
law then gives convergence of the full sequence and propagation
of chaos (Theorem~\ref{t-main}).

\subsection{Notation and tools}\label{ss-iset}
We introduce the canonical notation and two measure-theoretic
tools used to pass particle identities to the limit. The first
extends identities from cylinder functionals, and the second
passes identities to limits of random measures. Their proofs are
given in Appendix~\ref{a-measure-tools} for completeness.

Recall from Section~\ref{ss-paths} the canonical spaces
\begin{align*}
  \Om_1 = D([0,\infty),\R) \times C([0,\infty),\R)
  \times C([0,\infty),[0,\infty))
  \end{align*}
and
\begin{align*}
  \Om_2 = \Om_1 \times \Om_1 \times C_\uparrow([0,\infty)) ,
\end{align*}
where $C_\uparrow([0,\infty))$ is the set of continuous nondecreasing paths
starting at $0$. On $\Om_1$ we use the coordinate process
$(X, B, \ld)$ and write
$Y_s := (X_s, B_s, \ld_s)$. On $\Om_2$ we use the coordinate process
$(Y^1, Y^2, L)$, where 
$Y^1 = (X^1, B^1, \ld^1)$ and  
$Y^2 = (X^2, B^2, \ld^2)$.

Let $\sT \subseteq [0,\infty)$ be a countable dense set 
(to be specified further in Section~\ref{ss-ilim}). For $m \in \N$,
we fix a countable set
$\sC_m \subseteq C_c((\R \times \R \times [0,\infty))^m)$ that is
dense for the supremum norm. A \define{cylinder functional} is a
function $\Phi(Y) = \phi(Y_{s_1}, \dots, Y_{s_m})$ on $\Om_1$, with
$m \in \N$, $s_1 < \dots < s_m$ in $\sT$ and $\phi \in \sC_m$. We call
$s_1, \dots, s_m$ its \define{times}. There are countably many
cylinder functionals, each is bounded, and each is $J_1$-continuous at
every path whose first component does not jump at its times (see 
\cite[VI.2.3]{JacShi03}).

\begin{lemma}\label{l-ext}
  Let $Y$ be an $\Om_1$-valued random element, and let $\xi$ be an
  integrable random variable.
  If $s \in \sT \cup \set{\infty}$ and $\ee*{\Phi(Y) \, \xi} = 0$
for every cylinder functional $\Phi$ with times in $[0,s]$, then
$\ee*{\xi \giv Y_v, \ v \leq s} = 0$. 
\end{lemma}

The proof is given in Appendix~\ref{ap-cylinder-extension}.

A measurable function $F$ on a space $S$ is said to be \define{uniformly
integrable} with respect to a sequence $\seq[N]{m_N}$ of probability
measures on $S$ if $\sup_N \int (\abs{F} - K)^+ \, \diffd m_N \to 0$ as
$K \to \infty$.

\begin{lemma}\label{l-pass}
  Let $S$ be a Polish space, and let $\seq[N]{M_N}$ be a sequence
  of random probability measures on $S$ whose laws converge weakly to some
  $\mathfrak M \in \sP(\sP(S))$. Let $m_N(A) := \ee*{M_N(A)}$ and
  $m(A) := \int \rho(A) \, \mathfrak M(\diffd\rho)$ be the mean measures,
  and let $F \colon S \to \R$ be measurable
  and continuous at $m$-almost every point.
  \begin{enumerate}[(a)]
  \item\label{i-pass-pos} If $F \geq 0$, then
\begin{align*}
  \int \ee^{\rho}{F} \, \mathfrak M(\diffd\rho)
  \leq \liminf_N \int F \, \diffd m_N.
\end{align*}
 The same
    holds, without the continuity assumption, if 
    $F \colon S \to [0,\infty]$ is
    lower semicontinuous.
  \item\label{i-pass-zero} If $F$ is uniformly integrable with
    respect to $\seq[N]{m_N}$ and $\ee*{\abs*{\inner*{M_N}{F}}} \to 0$, then,
    for $\mathfrak M$-almost every $\rho$,
\begin{align*}
  \ee^{\rho}{\abs{F}} < \infty
\eand \ee^{\rho}{F} = 0.
\end{align*}
  \end{enumerate}
\end{lemma}

The proof is given in Appendix~\ref{a-random-measures}.

\subsection[The limit]{The limit}\label{ss-ilim}
By Proposition~\ref{p-tight}, the laws of a subsequence of $\seq[N]{R_N}$
converge weakly to some $\mathfrak R \in \sP(\sP(\Om_2))$. We fix such
a subsequence, as well as $\mathfrak R$, and proceed, notationally,
as if the full sequence converges. The empirical measures $R_N$ and
$Q_N$ are random probability measures on $\Om_2$ and $\Om_1$, so their
integrals are written with brackets: for bounded or nonnegative
measurable $F$ on $\Om_2$ and $G$ on $\Om_1$,
\begin{align*}
  \inner{R_N}{F} = \oo{N(N-1)} \sum_{i \neq j} F(\bar Y^{i,j,N}),
  \qquad
  \inner{Q_N}{G} = \oo{N} \sum_{i=1}^N G(Y^{i,N}).
\end{align*}
Deterministic elements of
$\sP(\Om_2)$ are denoted by $\rho$, and for such $\rho$ we put
$\rho_1 := \pushf{\pi_1}{\rho}$, the first marginal of $\rho$. If $F = G(Y^1)$ for a bounded or
nonnegative measurable $G$ on $\Om_1$, we write $\ee^{\rho_1}{F}$ for
$\ee^{\rho_1}{G}$, and $\inner{Q_N}{F}$ for $\inner{Q_N}{G}$. Then
$\ee^{\rho}{F} = \ee^{\rho_1}{F}$ and, since
$Q_N = \pushf{\pi_1}{R_N}$, also $\inner{R_N}{F} = \inner{Q_N}{F}$.

Let $\bar\rho(A) := \int \rho(A) \, \mathfrak R(\diffd\rho)$ be
the mean measure, or barycenter, of $\mathfrak R$. We apply Lemma~\ref{l-pass} with
$S = \Om_2$, $M_N = R_N$ and $\mathfrak M = \mathfrak R$. By
\eqref{mean}, $m_N$ is then the law of $\bar Y^{1,2,N}$, so that
$\int F \, \diffd m_N = \ee*{F(\bar Y^{1,2,N})}$, and $m = \bar\rho$. For
bounded continuous $F$ on $\Om_2$, the map
$\rho \mapsto \ee^{\rho}{F}$ is bounded and continuous. So weak
convergence and \eqref{mean} give
\begin{align*}
  \int F \, \diffd\bar\rho = \int \ee^{\rho}{F} \, \mathfrak R(\diffd\rho)
  = \lim_N \ee{\inner{R_N}{F}}
  = \lim_N \ee{F(\bar Y^{1,2,N})} ,
\end{align*}
and $\bar\rho$ is the weak limit of the laws of $\bar Y^{1,2,N}$.
So $\bar\rho$ is the limit of the annealed laws of a random
pair, while $\mathfrak R$ describes the limits of the quenched laws.

Let $\sJ_{\mathfrak R}$ be the countable set of $t > 0$ at which
$X^1$ or $X^2$ jumps with positive $\bar\rho$-probability. If
$t \notin \sJ_{\mathfrak R}$, then, for $\mathfrak R$-almost every
$\rho$, neither $X^1$ nor $X^2$ jumps at $t$, $\rho$-almost surely.
From now on, $\sT$ is a countable dense subset of
$[0,\infty) \setminus \sJ_{\mathfrak R}$ with $0 \in \sT$. Then, for
every cylinder functional $\Phi$, $\Phi(Y^1)$ and $\Phi(Y^2)$ are
continuous at $\bar\rho$-almost every point of $\Om_2$. Finally,
by Proposition~\ref{p-tight}\ref{i-tight-prod}, for
$\mathfrak R$-almost every $\rho$, the coordinates $Y^1$ and $Y^2$
are independent under $\rho$, each with law $\rho_1$.

\subsection{One particle}\label{ss-ione}

\begin{proposition}\label{p-ione}
  For $\mathfrak R$-almost every $\rho$, the law
  $\rho_1 = \pushf{\pi_1}{\rho}$ satisfies
  \ref{i-nl-init}--\ref{i-nl-mom} of Definition~\ref{d-nl}, and
  $\ee^{\rho_1}*{G(\sup_{s \leq T} \ld_s)} < \infty$ for every
  $T \geq 0$.
\end{proposition}

\begin{proof}
  Each step below holds outside an $\mathfrak R$-null set of
  $\rho$, and
  there are countably many steps. Under $\rho_1$, we use the coordinate
  process $(X, B, \ld)$ of $\Om_1$ and its jump measure
  $\mu := \Gamma(X)$.

  \emph{Initial law.} The map
  $\psi \colon \sP(\Om_2) \to \sP(\R \times [0,\infty))$,
  $\psi(\rho) := \pushf{(X_0, \ld_0)}{\rho_1}$, is continuous,
  since evaluation at time $0$ is continuous on $\Om_1$. The
  random measure $\psi(R_N) = \pushf{(X_0, \ld_0)}{Q_N}$ is the empirical measure of $(X^{i,N}_0, \ld^{i,N}_0)_{i \in [N]}$.
  Next, we use Assumption~\ref{a-std}\ref{i-eta}--\ref{i-chaos}
  as well as 
  \cite[Proposition 2.2]{Szn91} which states that for exchangeable particles
  chaos is equivalent to convergence in law of the empirical measures
  to a constant. Together, they imply that this empirical measure converges in
  law to the constant $\eta_0$. So $\psi(\rho) = \eta_0$ for
  $\mathfrak R$-almost every $\rho$, by the continuous mapping
  theorem. That is, $(X_0, \ld_0)$ has law $\eta_0$ under $\rho_1$,
  which is \ref{i-nl-init}. In the same way,
  $B_0 = 0$, $\rho_1$-almost surely.

  \emph{Moments and jumps.} The functional $G(\sup_{s \leq T} \ld_s)$
  is continuous on $\Om_1$. The number $n_T$ of jumps of $X$ in $(0,T)$
  is lower semicontinuous on $\Om_1$. Indeed, let $x^n \to x$ in
  $D([0,\infty),\R)$, and let $s_1 < \dots < s_k$ be jump times of $x$
  in $(0,T)$. By \cite[Theorem VI.1.14]{JacShi03}, there are
  continuous increasing bijections $\lambda_n$ of $[0,\infty)$ with
  $\sup_s \abs{\lambda_n(s) - s} \to 0$ and $x^n \circ \lambda_n \to x$
  uniformly on $[0,T]$. The jump of $x^n \circ \lambda_n$ at $s_i$ is
  the jump of $x^n$ at $\lambda_n(s_i)$, and it converges to
  $\Delta x(s_i) \neq 0$. So, for large $n$, $x^n$ jumps at the $k$
  distinct times $\lambda_n(s_1), \dots, \lambda_n(s_k)$, which lie in
  $(0,T)$, and $\liminf_n n_T(x^n) \geq k$. By
  Lemma~\ref{l-pass}\ref{i-pass-pos},
  Proposition~\ref{p-moments} and
  Lemma~\ref{l-jumpmom}\ref{i-bdnu} with $A = E$,
  \begin{align*}
    \int \ee^{\rho_1}{G\prn{\textstyle\sup_{s \leq T} \ld_s} + n_T} \,
    \mathfrak R(\diffd\rho) \leq C_T .
  \end{align*}
  This gives the $G$-moment, hence
  $\ee^{\rho_1}*{\sup_{s \leq T} \ld_s} \leq x_G
  + \ee^{\rho_1}*{G(\sup_{s \leq T} \ld_s)} < \infty$ and
  \ref{i-nl-mom}, and
  $\ee^{\rho_1}*{n_T} < \infty$ for all $T$. In particular, $X$ has
  finitely many jumps on bounded intervals, $\rho_1$-almost surely.

  \emph{The jumps of $X$.} Let $t \in \sT$ with $t > 0$, and
  $0 < \dl < M$. Let $h(z) := z \chi(z)$, where $\chi \in C_c(E)$
  takes values in $[0,1]$, equals $1$ on
  $\sets{z}{2\dl \leq \abs{z} \leq M}$, and vanishes on
  $\sets{z}{\abs{z} \leq \dl}$. Put
  \begin{align*}
    G := 1 \wedge \abs*{X_t - X_0 - B_t
    - \tsum_{0 < s \leq t} h(\Delta X_s)} .
  \end{align*}
  By Lemma~\ref{l-jump}\ref{i-jump-sum}, $G$ is continuous at every
  path whose first component does not jump at $t$, hence
  $\bar\rho$-almost everywhere. For the particle system, \eqref{X} gives
  $X^1_t - X^1_0 - B^1_t = \sum_{0 < s \leq t} \Delta X^1_s$, and
  $\abs{z - h(z)} \leq 2\dl$ for $\abs{z} \leq M$. So
  $G(Y^{1,N}) \leq 2\dl \, J^{1,N}_t$, unless $X^{1,N}$ has a jump of
  size larger than $M$ in $(0,t]$. By
  Lemma~\ref{l-jumpmom}\ref{i-bdnu},
  \begin{align*}
    \ee{G(Y^{1,N})}
    \leq 2\dl \, \ee*{J^{1,N}_t}
    + \ee{\mu^{1,N}\prn*{(0,t] \times \sets{z}{\abs{z} > M}}}
    \leq C_t \prn{\dl + \nu\prn*{\sets{z}{\abs{z} > M}}} .
  \end{align*}
  By Lemma~\ref{l-pass}\ref{i-pass-pos},
  $\int \ee^{\rho_1}{G} \, \mathfrak R(\diffd\rho)$ obeys
  the same bound. Let $\dl \downarrow 0$ and $M \uparrow \infty$ along
  sequences. Since $X$ has finitely many jumps on $(0,t]$,
  $\sum_{0 < s \leq t} h(\Delta X_s) \to
  \sum_{0 < s \leq t} \Delta X_s$, and Fatou's lemma gives
  $X_t - X_0 - B_t = \sum_{0 < s \leq t} \Delta X_s$,
  $\rho_1$-almost surely. By right continuity, this holds for all $t$ at
  once. So $X - X_0 - B$ is adapted and piecewise constant, starts at
  $0$, and has finitely many jumps on bounded intervals. Once
  \ref{i-nl-bm} is shown, this is \ref{i-nl-graph}.

  \emph{Brownian motion.} Let $s < t$ in $\sT$, $f \in C^2_c(\R)$, and
  let $\Phi$ be a cylinder functional with times in $[0,s]$. Put
  \begin{align*}
    F := \Phi(Y^1) \prn{f(B^1_t) - f(B^1_s)
    - \ot \int_s^t f''(B^1_v) \, \dv} .
  \end{align*}
  It is bounded and continuous $\bar\rho$-almost everywhere. For the
  particle system, the processes
  $f(B^i) - \ot \int_0^\cdot f''(B^i_v) \, \dv$ are martingales,
  bounded on bounded intervals, and orthogonal for distinct $i$, since
  $(B^1, \dots, B^N)$ is an $N$-dimensional $\bF$-Brownian motion. As
  $\Phi(Y^{i,N})$ is $\sF_s$-measurable, the cross terms in
  $\ee*{\inner*{Q_N}{F}^2}$ vanish, so
  $\ee*{\inner*{Q_N}{F}^2}$ $\leq \norm{F}_\infty^2 / N$. By
  Lemma~\ref{l-pass}\ref{i-pass-zero}, $\ee^{\rho_1}{F} = 0$. We let $f$
  range over a countable set that is dense in $C^2_c(\R)$ for the
  norm $\norm{f}_\infty + \norm{f'}_\infty + \norm{f''}_\infty$. By
  Lemma~\ref{l-ext}, by continuity in $s$ and $t$, and by
  approximation in $f$,
  \begin{align*}
    f(B) - \ot \int_0^\cdot f''(B_v) \, \dv
  \end{align*}
  is an $\bF^{\rho_1}$-martingale for every $f \in C^2_c(\R)$. Taking $f$
  equal to $x$ and to $x^2$ on $[-n,n]$, and stopping at the exit time
  of $(-n,n)$, we see that $B$ is a continuous local martingale with
  $B_0 = 0$ and $\qv{B}_t = t$. By L\'evy's characterization, $B$ is
  an $\bF^{\rho_1}$-Brownian motion. This is \ref{i-nl-bm}.

  \emph{Compensator.} Let $s < t$ in $\sT$, $h \in C_c(E)$, and let
  $\Phi$ be a cylinder functional with times in $[0,s]$. With
  $\nu(h) := \int_E h \, \diffd\nu$, put
  \begin{align*}
    F := \Phi(Y^1) \prn{\tsum_{s < v \leq t} h(\Delta X^1_v)
    - \nu(h) \int_s^t \ld^1_v \, \dv} .
  \end{align*}
  By Lemma~\ref{l-jump}\ref{i-jump-sum}, $F$ is continuous
  $\bar\rho$-almost everywhere. Since
  $\abs{F} \leq \norm{\Phi}_\infty \norm{h}_\infty
  (J^1_t + \int_0^t \ld^1_v \, \dv)$ for the particle system,
  Lemma~\ref{l-jumpmom}\ref{i-bdjsq}
  shows that the variables $F(Y^{1,N})$, $N \geq 2$, are
  uniformly integrable. By
  condition~\ref{i-comp} of Definition~\ref{d-sol}, the processes
  \begin{align*}
    M^i := \tsum_{0 < v \leq \cdot} h(\Delta X^i_v)
    - \nu(h) \int_0^\cdot \ld^i_v \, \dv
  \end{align*}
  are martingales with
  $\ee*{(M^i_t)^2} = \ee*{\int_0^t \ld^i_v \, \dv \int h^2 \, \diffd\nu}
  \leq C_t$. For $i \neq j$, $M^i$ and $M^j$ are purely discontinuous
  and have no common jumps, by condition~\ref{i-nosim}. So $M^i M^j$
  is a local martingale. On $[0,t]$ it is dominated by
  $\sup_{v \leq t} \abs{M^i_v} \sup_{v \leq t} \abs{M^j_v}$, which is
  integrable by Doob's inequality, so it is a martingale. As before,
  the cross terms vanish, and
  $\ee*{\inner*{Q_N}{F}^2}
  \leq \norm{\Phi}_\infty^2 \, \ee*{(M^1_t - M^1_s)^2}/N
  \leq C_t/N \to 0$. By
  Lemma~\ref{l-pass}\ref{i-pass-zero}, $\ee^{\rho_1}{F} = 0$.

  We let $h$ range over a countable set $\sC_E \subseteq C_c(E)$ such
  that every $h \in C_c(E)$ is a uniform limit of elements of $\sC_E$
  supported in a fixed compact subset of $E$. Since
  $\ee^{\rho_1}*{n_{t+1} + \int_0^t \ld_v \, \dv} < \infty$,
  Lemma~\ref{l-ext}, right continuity in $s$ and $t$, and approximation in $h$
  show that, for all $s < t$, $C \in \sF^{\rho_1}_s$ and $h \in C_c(E)$,
  \begin{align*}
    \ee^{\rho_1}{\ind{C} \int_{(s,t] \times E} h \, \diffd\mu}
    = \nu(h) \, \ee^{\rho_1}{\ind{C} \int_s^t \ld_v \, \dv} .
  \end{align*}
  Both sides are integrals of $h$ against finite measures on $E$, so
  the identity extends to $h = \ind{A}$ for every Borel
  $A \subseteq E$. The sets $C \times (s,t] \times A$
  generate the predictable $\ssig$ on
  $\Om \times (0,\infty) \times E$, and on them the measures induced
  by $\mu$ and by $\ld_s \, \ds \, \nu(\dz)$ agree. Hence
  $\ld_s \, \ds \, \nu(\dz)$ is the $\bF^{\rho_1}$-compensator of $\mu$
  \cite[Theorem II.1.8]{JacShi03}. This is \ref{i-nl-comp}.
\end{proof}

\subsection{The collision local time}\label{ss-ipair}
Suppose that $\rho_1$ satisfies \ref{i-nl-bm} and \ref{i-nl-graph}, and
that $(Y^1, Y^2)$ has law $\rho_1 \otimes \rho_1$. In the usual augmentation of
the filtration generated by $Y^1$ and $Y^2$, the process
$Z := X^1 - X^2$ is then piecewise Brownian with variance $2$, with
$W := B^1 - B^2$. We write $\ell(Z)$ for its local time at
$0$, defined by \eqref{tanaka} with $c = 2$. To have a Borel function
on all of $\Om_2$, we use the version
$\ell_t(Z) := \limsup_n \frac{n}{2} \int_0^t
\inds{\abs*{X^1_s - X^2_s} < 1/n} \, \ds$, which agrees with the
local time almost surely, by Section~\ref{ss-loctime}.

\begin{proposition}\label{p-ipair}
  For $\mathfrak R$-almost every $\rho$, under $\rho$ the pair $(Y^1, Y^2)$ has law
  $\rho_1 \otimes \rho_1$, and $L = \ell(Z)$ almost surely.
\end{proposition}

\begin{proof}
  The first claim is Proposition~\ref{p-tight}\ref{i-tight-prod}. By
  Proposition~\ref{p-ione}, $\ell(Z)$ is then defined under
  $\rho$, for $\mathfrak R$-almost every $\rho$. Fix $t \in \sT$ and $\eps \in (0,1]$. Let
  $k_\eps(a)$ $:= \eps^{-1} (1 - \abs{a}/\eps)^+$, and put
  \begin{align*}
    F_\eps := 1 \wedge \abs*{L_t
    - \int_0^t k_\eps(X^1_s - X^2_s) \, \ds} .
  \end{align*}
  Convergence in $\Om_1 \times \Om_1$ implies convergence of the first
  components at all but countably many times. So $F_\eps$ is
  continuous on $\Om_2$, by bounded convergence. For the particle
  system, let $\ell_t(a)$ be the local time of $Z^{1,2,N}$ at
  $a$, so that $\ell_t(0) = \ell^{1,2,N}_t$ and
  $\int_0^t k_\eps(Z^{1,2,N}_s) \, \ds
  = \int \ell_t(a) \, k_\eps(a) \, \da$. The number of pieces of
  $Z^{1,2,N}$ on $[0,t)$ is at most $1 + J^{1,N}_t + J^{2,N}_t$, which
  has expectation at most $C_t$. By
  Lemma~\ref{l-pieces}\ref{i-pieces-lev}, with $\sg = 0$,
  $\tau = \th = t$ and $c = 2$,
  \begin{align*}
    \ee{F_\eps(\bar Y^{1,2,N})}
    \leq \int \ee{\abs*{\ell_t(a) - \ell_t(0)}} \, k_\eps(a) \, \da
    \leq C_t \, \eps^{1/2} .
  \end{align*}
  By Lemma~\ref{l-pass}\ref{i-pass-pos},
  $\int \ee^{\rho}{F_\eps} \, \mathfrak R(\diffd\rho)$
  $\leq C_t \, \eps^{1/2}$. Under $\rho$, the
  local time of $Z$ is continuous in the level, so
  $\int_0^t k_\eps(Z_s) \, \ds \to \ell_t(Z)$ as $\eps \to 0$. By
  Fatou's lemma, applied along $\eps = 1/n$,
  $\int \ee^{\rho}{1 \wedge \abs{L_t - \ell_t(Z)}} \,
  \mathfrak R(\diffd\rho) = 0$. So
  $L_t = \ell_t(Z)$, $\rho$-almost surely, for
  $\mathfrak R$-almost every $\rho$, for all $t \in \sT$, and by
  continuity for all $t$.
\end{proof}

\subsection{The intensity}\label{ss-ifeed}
For $l \in C_\uparrow([0,\infty))$ and $t \geq 0$, put
\begin{align*}
  \Theta_t(l) := \int_0^t e^{-(t-s)} \, \diffd l_s
  = l_t - \int_0^t e^{-(t-s)} \, l_s \, \ds .
\end{align*}
The map $\Theta_t$ is continuous on $C_\uparrow([0,\infty))$, and
$0 \leq \Theta_t(l) \leq l_t$. In this notation, \eqref{ld} reads
$\ld^i_t = e^{-t} \ld^i_0 + \gm \sum_{j \neq i} \Theta_t(\ell^{i,j})$.

\begin{theorem}\label{t-ident}
  For $\mathfrak R$-almost every $\rho$, the law
  $\rho_1 = \pushf{\pi_1}{\rho}$ is a solution of the nonlinear system, and
  $\rho$ is the law of $(Y^1, Y^2, \ell(X^1 - X^2))$, where $Y^1$ and
  $Y^2$ are independent with law $\rho_1$.
\end{theorem}

\begin{proof}
  By Propositions~\ref{p-ione} and~\ref{p-ipair}, it remains to prove
  \ref{i-nl-feed} for $\mathfrak R$-almost every $\rho$. By Proposition~\ref{p-ione}
  and Lemma~\ref{l-dens}, $X_t$ has a continuous density $u(t,\cdot)$
  under $\rho_1$, for every $t > 0$.

  Fix $t \in \sT$ and a cylinder functional $\Phi$, and put
  \begin{align*}
    H := \Phi(Y^1) \prn{\ld^1_t - e^{-t} \ld^1_0 - \kp \Theta_t(L)} .
  \end{align*}
  It is continuous $\bar\rho$-almost everywhere, and
  $\abs{H} \leq \norm{\Phi}_\infty (\ld^1_t + \ld^1_0 + \kp L_t)$. By
  Corollaries~\ref{c-mom}\ref{i-mom-ui} and~\ref{c-ltui}, the
  variables $H(\bar Y^{1,2,N})$ are uniformly integrable. For the particle
  system,
  \begin{align*}
    \inner{R_N}{H} = \frac{(N-1)\gm_N - \kp}{N(N-1)}
    \sum_{i \neq j} \Phi(Y^{i,N}) \, \Theta_t(\ell^{i,j,N}),
  \end{align*}
  so, by \eqref{mean},
  $\ee*{\abs*{\inner*{R_N}{H}}}$ $\leq \norm{\Phi}_\infty
  \abs*{(N-1)\gm_N - \kp} \, \ee*{\ell^{1,2,N}_t}$. This tends to $0$,
  since $(N-1)\gm_N \to \kp$ by Assumption~\ref{a-std}\ref{i-mf}, and
  $\ee*{\ell^{1,2,N}_t} \leq C_t$ by Proposition~\ref{p-aldous},
  applied on intervals of length at most $1$. By
  Lemma~\ref{l-pass}\ref{i-pass-zero}, $\ee^{\rho}{H} = 0$ for
  $\mathfrak R$-almost every $\rho$.

  Fix $\rho$ outside the exceptional sets, so that
  $\ee^{\rho}{H} = 0$
  for all $t \in \sT$ and all cylinder functionals $\Phi$. By
  Proposition~\ref{p-ipair}, $L = \ell(Z)$ under $\rho$. By
  Lemma~\ref{l-pass}\ref{i-pass-pos}, $\ld^1_t$, $\ld^1_0$ and $L_t$
  are $\rho$-integrable. Lemma~\ref{l-ext}, with $s = \infty$, gives
  \begin{align*}
    \ld^1_t - e^{-t} \ld^1_0
    = \kp \, \ee^{\rho}{\Theta_t(\ell(Z)) \giv Y^1}
    \quad \rho\text{-a.s.}
  \end{align*}
  By Proposition~\ref{p-ione}, Lemma~\ref{l-condlt} applies to $\rho_1$,
  so $\ee^{\rho}{\ell_s(Z) \giv Y^1}$ $= \int_0^s u(v, X^1_v) \, \dv$ for
  every $s \geq 0$. As $\ell(Z)$ is continuous, nondecreasing and
  integrable, Fubini's theorem for conditional expectations and an
  integration by parts give
  \begin{align*}
    \ee^{\rho}{\Theta_t(\ell(Z)) \giv Y^1}
    &= \int_0^t u(s, X^1_s) \, \ds
    - \int_0^t e^{-(t-s)} \int_0^s u(v, X^1_v) \, \dv \, \ds
 \\ &= \int_0^t e^{-(t-s)} u(s, X^1_s) \, \ds .
  \end{align*}
  So \eqref{feed} holds $\rho_1$-almost surely for every $t \in \sT$. Both
  sides are continuous in $t$, pathwise, since
  $\int_0^t \norm{u(s,\cdot)}_\infty \, \ds < \infty$ by the bound
  in the proof of Lemma~\ref{l-condlt}. Hence \eqref{feed} holds for
  all $t$ at once.
\end{proof}

\subsection{Propagation of chaos}\label{ss-chaos}

We now combine tightness, the identification of the limit and
uniqueness. Recall the Borel version $\ell(X^1 - X^2)$ of the
collision local time on $\Om_2$ from Section~\ref{ss-ipair}.

\begin{theorem}\label{t-main}
  There is exactly one solution $P$ of the nonlinear system. Let
  $P_2 \in \sP(\Om_2)$ be the law of $(Y^1, Y^2, \ell(X^1 - X^2))$,
  where $Y^1$ and $Y^2$ are independent with law $P$. Then, as
  $N \to \infty$:
  \begin{enumerate}[(a)]
  \item\label{i-main-prob} $Q_N \to P$ in probability in
    $\sP(\Om_1)$, and $R_N \to P_2$ in probability in $\sP(\Om_2)$;
  \item\label{i-main-chaos} for every $k \in \N$, the law of
    $(Y^{1,N}, \dots, Y^{k,N})$ converges weakly to $P^{\otimes k}$;
  \item\label{i-main-pair} the law of $\bar Y^{1,2,N}$ converges
    weakly to $P_2$.
  \end{enumerate}
\end{theorem}

\begin{proof}
  \ref{i-main-prob} By Proposition~\ref{p-tight}\ref{i-tight-R} and
  Prokhorov's theorem, every subsequence of $\seq[N]{R_N}$ has a
  further subsequence along which the laws of $R_N$ converge weakly
  to some $\mathfrak R \in \sP(\sP(\Om_2))$. Section~\ref{s-ident}
  applies to every such $\mathfrak R$. By Theorem~\ref{t-ident}, for
  $\mathfrak R$-almost every $\rho$, the first marginal $\rho_1$ is a
  solution of the nonlinear system, and $\rho$ is the law of
  $(Y^1, Y^2, \ell(X^1 - X^2))$ with $Y^1$ and $Y^2$ independent with
  law $\rho_1$. In particular, a solution exists. By
  Proposition~\ref{p-nlunique}, it is unique; call it $P$. Hence
  $\rho_1 = P$ and $\rho = P_2$ for $\mathfrak R$-almost every
  $\rho$, that is, $\mathfrak R = \delta_{P_2}$. As the limit does not
  depend on the subsequence, the laws of $R_N$ converge weakly to
  $\delta_{P_2}$. Convergence in law to a constant is convergence in
  probability, so $R_N \to P_2$ in probability. Since
  $Q_N = \pushf{\pi_1}{R_N}$, $\pushf{\pi_1}{P_2} = P$, and the map
  $\rho \mapsto \pushf{\pi_1}{\rho}$ is continuous, also
  $Q_N \to P$ in probability.

  \ref{i-main-chaos} By exchangeability of the finite particle system, the random elements
  $Y^{1,N}, \dots, Y^{N,N}$ of $\Om_1$ are also exchangeable. By
  \cite[Proposition 2.2]{Szn91}, convergence of their empirical
  measures $Q_N$ in law to the constant $P$, which is part of
  \ref{i-main-prob}, is equivalent to \ref{i-main-chaos}.

  \ref{i-main-pair} Let $F$ be bounded and continuous on $\Om_2$.
  The map $\rho \mapsto \ee^{\rho}{F}$ is bounded and continuous on
  $\sP(\Om_2)$. So \eqref{mean} and \ref{i-main-prob} give
  $\ee*{F(\bar Y^{1,2,N})} = \ee*{\inner{R_N}{F}} \to \ee^{P_2}{F}$.
\end{proof}

The limit keeps genuine jumps, while their fluctuations average out
in the empirical measures.

\begin{remark}\label{r-jumps}
  Suppose that $\kp > 0$. Under $P$, $X_s$ has the density
  $u(s,\cdot)$ for $s > 0$, so
  $P\prn*{u(s, X_s) = 0} = \int_{\set*{u(s,\cdot) = 0}} u(s,x) \, \dx = 0$.
  By Fubini's theorem, $P$-almost surely $u(s, X_s) > 0$ for almost
  every $s$, and then \eqref{feed} gives $\ld_t > 0$ for all $t > 0$.
  Choose $0 < a < b$ such that $G := \sets*{z}{a < \abs{z} \leq b}$
  has $\nu(G) > 0$, and put
  \begin{align*}
    M^G_t := \int_{(0,t] \times G} z \prn*{\mu(\ds,\dz)
    - \ld_s \, \ds \, \nu(\dz)} .
  \end{align*}
  Such $a$ and $b$ exist since $\nu(E) = 1$, and restricting to
  jumps with sizes in $G$ avoids any moment assumption on $\nu$.
  By \ref{i-nl-comp} and \ref{i-nl-mom} of Definition~\ref{d-nl},
  $M^G$ is a square-integrable martingale under $P$ with
  \begin{align*}
    \qv{M^G}_t = \int_G z^2 \, \nu(\dz) \int_0^t \ld_s \, \ds > 0
    \quad \text{for } t > 0, \ P\text{-a.s.}
  \end{align*}
  So the jumps of $X$ do not disappear in the limit. For the particle
  system, let $M^{G,i,N}$ be the same martingale built from $\mu^{i,N}$
  and $\ld^{i,N}$, and $M^{G,N} := \oo{N} \sum_{i=1}^N M^{G,i,N}$. Two
  particles never jump at the same time, so the $M^{G,i,N}$ are
  orthogonal, and Corollary~\ref{c-mom}\ref{i-mom-sup} gives
  \begin{align*}
    \ee{(M^{G,N}_t)^2} = \oo{N^2} \int_G z^2 \, \nu(\dz)
    \sum_{i=1}^N \ee{\int_0^t \ld^{i,N}_s \, \ds}
    \leq \frac{b^2 C_t}{N} .
  \end{align*}
  So the compensated jumps vanish after averaging over the particles,
  in line with the deterministic limit of $Q_N$.
\end{remark}

\appendix
\section{Construction and Poisson representations}\label{a-construction}
We give the construction and uniqueness proofs deferred from
Section~\ref{s-model}, using the notation of that section.

\subsection{Construction up to explosion}\label{ap-exist}
\begin{proof}[Proof of Proposition~\ref{p-exist}]
  Between jumps, the positions are known Brownian paths shifted by
  constants. So the pair local times and the intensities are
  determined, and the next jump is the first atom of some $\Pi^i$
  below the graph of $\ld^i$. We only indicate the main steps.

  Since $\nu(E) = 1$, the atoms of $\Pi^1, \dots, \Pi^N$ with level at
  most $M$ have time coordinates forming a Poisson process of rate
  $NM$. Discarding a null set, we may therefore assume that, for every
  $M \in \N$, each $\Pi^i$ has finitely many atoms in
  $(0,M] \times E \times (0,M]$, and that no two atoms share a time
  coordinate.

  Put $T_0 := 0$. Suppose that the first $n$ jump times
  $T_1, \dots, T_n$ have been constructed, together with the indices
  $\iota_1, \dots, \iota_n \in [N]$ of the jumping particles and the
  jump sizes $\xi_1, \dots, \xi_n \in E$. For $t \geq 0$, let
  \begin{align*}
    X^{n,i}_t := X^i_0 + B^i_t
    + \sum_{k=1}^n \xi_k \inds{\iota_k = i,\ T_k \leq t},
  \end{align*}
  and let $\ld^{n,i}$ be given by \eqref{ld}, with $\ell^{i,j}$
  the local time at $0$ of the piecewise Brownian
  process $X^{n,i} - X^{n,j}$, with variance $2$. The
  processes $\ld^{n,i}$ are continuous and adapted. Let $T_{n+1}$ be
  the first time after $T_n$ at which some $\Pi^i$ has an atom
  $(s,z,y)$ with $y \leq \ld^{n,i}_s$, and let $\iota_{n+1} := i$ and
  $\xi_{n+1} := z$. Since the $\ld^{n,i}$ are locally bounded, the
  choice of the null set gives $T_{n+1} > T_n$ on $\set*{T_n < \infty}$
  and uniqueness of this atom. Since the thinning integrand
  $\inds*{s > T_n,\ y \leq \ld^{n,i}_s}$ is predictable, $T_{n+1}$ is a
  stopping time and $(\iota_{n+1}, \xi_{n+1})$ is
  $\sF_{T_{n+1}}$-measurable.

  Local times depend only on the past and are continuous, and
  $X^{n+1} = X^n$ on $\llbracket 0, T_{n+1} \llbracket$. Hence
  $\ld^{n+1} = \ld^n$ on $\llbracket 0, T_{n+1} \rrbracket$. So
  $\zeta := \sup_n T_n$, $X^i := X^{n,i}$ and $\ld^i := \ld^{n,i}$ on
  $\llbracket T_n, T_{n+1} \llbracket$, and
  $\mu^i := \sum_{n \geq 1,\, T_n < \infty,\, \iota_n = i}
  \delta_{(T_n, \xi_n)}$ are
  well defined. Items \ref{i-space}--\ref{i-proc} and conditions \ref{i-nosim}
  and~\ref{i-eqs} of
  Definition~\ref{d-sol} hold by construction.

  On $\rrbracket T_n, T_{n+1} \rrbracket$ we have $\ld^i = \ld^{n,i}$,
  and the only atom of $\Pi^1, \dots, \Pi^N$ there below the graphs
  of the $\ld^i$ is the one at time $T_{n+1}$. This gives
  \eqref{thin}. Condition~\ref{i-comp} of Definition~\ref{d-sol} follows: for
  nonnegative predictable $H$, the function
  $H(s,z) \inds*{s \leq T_n,\ y \leq \ld^i_s}$ is predictable, so the
  compensator of $\Pi^i$ gives
  \begin{align*}
    \ee{\int H(s,z) \inds{s \leq T_n} \, \mu^i(\ds,\dz)}
    = \ee{\int_0^{T_n} \int_E H(s,z) \, \ld^i_s \, \nu(\dz) \, \ds}.
  \end{align*}
  Finally, suppose that $\zeta < \infty$ and
  $\sup_{t < \zeta} \sum_i \ld^i_t \leq M$ for an integer
  $M \geq \zeta$. By \eqref{thin}, the infinitely many jump times are
  time coordinates of atoms in $(0,M] \times E \times (0,M]$. This
  contradicts the choice of the null set and proves part~\ref{i-blowup} of the
  proposition.
\end{proof}

\subsection{Poisson representation and uniqueness in law}\label{ap-unique}
\begin{proof}[Proof of Proposition~\ref{p-unique}]
  By Proposition~\ref{p-nonexpl}, $\zeta = \infty$. For $T_n < \infty$,
  let $\iota_n \in [N]$ and $\xi_n \in E$ be the index of the particle
  that jumps at $T_n$ and its jump size. By condition~\ref{i-comp},
  $\ee*{\int \inds*{\ld^i_s = 0} \, \mu^i(\ds,\dz)} = 0$, so
  $\ld^{\iota_n}_{T_n} > 0$ for all $n$, a.s.

  Enlarge the space by random variables $U_1, U_2, \dots$, uniform on
  $(0,1]$, and Poisson random measures
  $\tilde\Pi^1, \dots, \tilde\Pi^N$ with intensity
  $\ds \, \nu(\dz) \, \dy$, all independent of each other and of $\sF$.
  Let $\bG$ be the usual augmentation of the filtration generated by
  $\sF_t$, the variables $U_n \inds*{T_n \leq t}$, and the restrictions
  of the $\tilde\Pi^i$ to $(0,t] \times E \times (0,\infty)$. The added
  objects are independent of $\sF$, and $U_n$ is revealed at the
  $\bF$-stopping time $T_n$. Hence every $\bF$-martingale is a
  $\bG$-martingale. So $B$ is a $\bG$-Brownian motion, condition~\ref{i-comp}
  holds with $\bG$ in place of $\bF$, and each $\tilde\Pi^i$ is a
  $\bG$-Poisson random measure. Set
  \begin{align*}
    \Pi^i(\ds,\dz,\dy) := \sum_{n:\, T_n < \infty,\, \iota_n = i}
    \delta_{(T_n, \xi_n, U_n \ld^i_{T_n})}
    + \inds{y > \ld^i_s} \, \tilde\Pi^i(\ds,\dz,\dy).
  \end{align*}
  Given $\sG_{T_n-}$ and $(\iota_n, \xi_n)$, the variable $U_n$ is
  still uniform on $(0,1]$. So, by condition~\ref{i-comp}, the
  $\bG$-compensator of the first term is
  $\inds*{y \leq \ld^i_s} \, \ds \, \nu(\dz) \, \dy$. The second term has
  compensator $\inds*{y > \ld^i_s} \, \ds \, \nu(\dz) \, \dy$, since
  $\ld^i$ is continuous and adapted. Almost surely, no two atoms of
  $\Pi^1, \dots, \Pi^N$ share a time coordinate. Regarded as one point
  process on $[N] \times E \times (0,\infty)$, the family
  $(\Pi^i)_{i \in [N]}$ therefore has a deterministic compensator. By
  the martingale characterization of Brownian motion and Poisson
  point processes \cite[Chapter II, Theorem 6.3]{IkeWat89},
  $B, \Pi^1, \dots, \Pi^N$ are a Brownian motion and Poisson random
  measures, independent of each other and of $\sG_0 \supseteq \sF_0$.

  The atoms of the first term in $\Pi^i$ lie below the graph of
  $\ld^i$, and those of the second term lie above it. This gives
  \eqref{thin}. By Remark~\ref{r-thin-unique}, $(X, \ld, \mu)$ is the
  solution of Proposition~\ref{p-exist} built from
  $(X_0, \ld_0, B, \Pi^1, \dots, \Pi^N)$. That solution is a
  measurable function of these drivers, which depends only on $\gm$.
  The law of the drivers depends only on $\eta$ and $\nu$.
\end{proof}

\section{An explicit uniform-integrability function}\label{a-ui-function}
We construct the functions $g$ and $G$ used in
Section~\ref{ss-bset} directly from the initial intensity laws.
Write $Y_N:=y_1$ under $\eta_N$, and set
\begin{align*}
  m_1&:=\sup_{N\geq2}\int Y_N\,\diffd\eta_N<\infty,\\
  H(r)&:=\sup_{N\geq2}\int Y_N\inds{Y_N>r}\,\diffd\eta_N.
\end{align*}
Uniform integrability gives $H(r)\to0$. For $n\geq1$ define
\begin{align*}
  r_n&:=\min\sets{r\in\N,\ r\geq1}{H(r)\leq2^{-n-1}},
  &a_n&:=2^n(1+m_1)r_n,
\end{align*}
and put, for $x\geq0$,
\begin{equation}\label{g-explicit}
  g(x):=\sum_{n=1}^{\infty}\frac{x}{a_n+x},
  \qquad
  G(x):=\sum_{n=1}^{\infty}
  \left(x-a_n\log\left(1+\frac{x}{a_n}\right)\right).
\end{equation}
These are explicit in terms of the common tail function $H$.

Since $a_n\geq2^n$, the series for $g$ and its first derivative
converge uniformly on compact intervals: for $0\leq x\leq R$,
\begin{align*}
  0\leq\frac{x}{a_n+x}\leq\frac{R}{a_n},\qquad
  0<\frac{a_n}{(a_n+x)^2}\leq\frac1{a_n}.
\end{align*}
Each summand is increasing and concave. Thus $g$ is nonnegative,
nondecreasing, concave and continuously differentiable on
$[0,\infty)$, with $g(0)=0$. For every integer $k\geq1$, the first
$k$ summands tend to $1$ as $x\to\infty$, so
$\liminf_{x\to\infty}g(x)\geq k$. Hence $g(x)\to\infty$.
Integrating term by term gives $G(x)=\int_0^x g(u)\,\diffd u$ and
the second formula in \eqref{g-explicit}. In particular,
$G(x)/x\geq g(x/2)/2\to\infty$.

It remains to verify the normalization \eqref{Gbound}. Let
$G_n(x):=x-a_n\log(1+x/a_n)$. Its integral representation gives
\begin{align*}
  0\leq G_n(x)=\int_0^x\frac{u}{a_n+u}\,\diffd u
  \leq\min\left\{x,\frac{x^2}{2a_n}\right\}.
\end{align*}
Splitting at $r_n$, for every $N\geq2$ we obtain
\begin{align*}
  \int G_n(Y_N)\,\diffd\eta_N
  &\leq\frac{r_n}{2a_n}\int Y_N\,\diffd\eta_N+H(r_n)\\
  &\leq\frac{m_1}{2^{n+1}(1+m_1)}+2^{-n-1}
  \leq2^{-n}.
\end{align*}
Tonelli's theorem now yields
\begin{align*}
  \sup_{N\geq2}\int G(y_1)\,\diffd\eta_N
  \leq\sum_{n=1}^{\infty}2^{-n}=1.
\end{align*}
This proves every property of $g$ and $G$ invoked in
Section~\ref{ss-bset}, without any second-moment assumption.

\section{Local-time estimates}\label{a-localtimes}
The following lemma applies to the piecewise Brownian processes
of Section~\ref{ss-loctime}. It controls local-time increments
through the number of Brownian pieces, without assumptions on
the sizes of the jumps.

\begin{lemma}\label{l-pieces}
  Let $Z = Z_0 + W + V$ be a piecewise Brownian process with
  variance $c$, and $\theta \geq 0$.

  Given two finite 
  $\bF$-stopping times $\sg,\tau$ such that 
  $\sg \leq \tau \leq \sg + \theta$, we set 
  $\sg_0 := \sg$, let $\sg_1 < \sg_2 < \cdots$ be the jump
  times of $V$ after $\sg$, and define 
\begin{align}
  \label{num-pieces}
  K := \tsum_{k \geq 0} \inds*{\sg_k < \tau}.
\end{align}
Then
  \begin{enumerate}[(a)]
  \item\label{i-pieces-bd} There are random variables
    $\zeta_k(a) \geq 0$, measurable in $(a, \om)$ and vanishing on
    $\set*{\sg_k = \infty}$, such that, for all $a \in \R$,
    \begin{align*}
      \ell_{\tau}(a) - \ell_{\sg}(a)
      \leq \sum_{k \geq 0} \inds{\sg_k < \tau} \, \zeta_k(a) \eand
      \ee{\zeta_k(a)^2 \giv \sF_{\sg_k}} \leq \theta / c .
    \end{align*}
    In particular, by Jensen's inequality,
    $\ee*{\ell_{\tau}(a) - \ell_{\sg}(a)} \leq \sqrt{\theta / c} \, \ee*{K}$.
  \item\label{i-pieces-ui} Fix $\theta$, $c$ and a function
    $\vp \colon \N \to [0,\infty)$ with $\vp(M) \to 0$ as
    $M \to \infty$. The random
    variables $\int (\ell_{\tau}(a) - \ell_{\sg}(a)) \, \rho(\da)$, over all
    $Z$, $\sg$, $\tau$ as above, on any filtered probability space, with
    $\ee*{(K - M)^+} \leq \vp(M)$ for all $M \in \N$, and all probability
    measures $\rho$ on $\R$, form a uniformly integrable family.
  \item\label{i-pieces-lev} For every $a \in \R$,
    \begin{align*}
      \ee{\abs{\prn{\ell_{\tau}(a) - \ell_{\sg}(a)} - \prn{\ell_{\tau}(0) - \ell_{\sg}(0)}}}
      \leq \prn{2 \abs{a} + 4 \abs{a}^{1/2} (c\th)^{1/4}}/c \,
      \, \ee*{K} .
    \end{align*}
  \end{enumerate}
\end{lemma}
\begin{proof}
  \ref{i-pieces-bd} On $\set*{\sg_k < \infty}$ put
  $\beta^k_r := Z_{\sg_k} + W_{\sg_k + r} - W_{\sg_k}$, $r \geq 0$.
  On $[\sg_k, \sg_{k+1})$, $Z$ follows $\beta^k$, and
  $\sg_{k+1} \wedge \tau - \sg_k \leq \theta$ on $\set*{\sg_k < \tau}$.
  Let $\zeta_k(a)$ be the local time of $\beta^k$ at $a$ on
  $[0,\theta]$. Since $t \mapsto \ell_t(a)$ is continuous, the bound
  on $\ell_{\tau}(a) - \ell_{\sg}(a)$ follows. Given $\sF_{\sg_k}$, $\beta^k$
  is $Z_{\sg_k}$ plus a Brownian motion with variance $c$. Its
  local time at $a$ on $[0,\theta]$ is
  stochastically largest for the start at $a$. There it is $1/c$
  times the semimartingale local time at $a$, which has the law of
  $\abs*{W_\theta}$ by L\'evy's theorem. So the conditional second
  moment is at most $\ee*{W_\theta^2}/c^2 = \theta/c$. The last claim
  follows since $\set*{\sg_k < \tau} \in \sF_{\sg_k}$ and
  $\ee*{\zeta_k(a) \giv \sF_{\sg_k}} \leq \sqrt{\theta / c}$.

  \ref{i-pieces-ui} Let $\xi := \int (\ell_{\tau}(a) - \ell_{\sg}(a)) \,
  \rho(\da)$ and $g_k := \int \zeta_k(a) \, \rho(\da)$. By
  \ref{i-pieces-bd}, Jensen's inequality and Fubini's theorem,
  $\xi \leq \sum_{k \geq 0} \inds*{\sg_k < \tau} \, g_k$ and
  $\ee*{g_k^2 \giv \sF_{\sg_k}} \leq \theta / c$, hence
  $\ee*{g_k \giv \sF_{\sg_k}} \leq \sqrt{\theta / c}$. For $M \in \N$,
  \begin{align*}
    \xi \leq \sum_{k < M} g_k
    + \sum_{k \geq M} \inds{\sg_k < \tau} \, g_k .
  \end{align*}
  The first sum has $\bL^2$ norm at most $M \sqrt{\theta / c}$. The second
  has expectation at most
  $\sqrt{\theta / c} \, \ee*{(K - M)^+} \leq \sqrt{\theta / c} \, \vp(M)$.
  Write these two sums as $A_M$ and $B_M$. For $R>0$,
  $(\xi-R)^+\leq(A_M-R)^++B_M$, so
  \begin{align*}
    \ee{(\xi-R)^+}
    \leq \frac{M^2\theta}{cR}+\sqrt{\theta/c}\,\vp(M).
  \end{align*}
  First let $R\to\infty$, then $M\to\infty$. This proves uniform
  integrability with only a common integrable-tail bound on $K$.

  \ref{i-pieces-lev} With $\beta^k$ as in \ref{i-pieces-bd}, let
  $\ell^k_r(a)$ be the local time of $\beta^k$ at $a$ on
  $[0,r]$, and
  $\xi_k := \sup_{r \leq \th} \abs*{\ell^k_r(a) - \ell^k_r(0)}$ on
  $\set*{\sg_k < \infty}$. As in \ref{i-pieces-bd}, the left side is
  at most $\ee*{\sum_{k \geq 0} \inds*{\sg_k < \tau} \, \xi_k}$. Given
  $\sF_{\sg_k}$, $\beta^k$ is $Z_{\sg_k}$ plus the Brownian motion
  $W^k := W_{\sg_k + \cdot} - W_{\sg_k}$ with variance $c$, which is
  independent of $\sF_{\sg_k}$. By
  \eqref{tanaka} at the levels $a$ and $0$,
  $c\,(\ell^k_r(a) - \ell^k_r(0))$ is the sum of
  $\prn{\abs*{\beta^k_r - a} - \abs*{\beta^k_r}}
  - \prn{\abs*{\beta^k_0 - a} - \abs*{\beta^k_0}}$, which is at most
  $2\abs{a}$ in absolute value, and of
  $-\int_0^r g_v \, \diffd W^k_v$ with
  $g_v := \sgn(\beta^k_v - a) - \sgn(\beta^k_v)$. The integrand is
  bounded by $2$ and vanishes unless $\beta^k_v$ lies between $0$ and
  $a$, which has conditional probability at most
  $\abs{a} (2\pi c v)^{-1/2}$. By Doob's inequality,
  \begin{align*}
    \ee{\sup_{r \leq \th} \abs*{\tint_0^r g_v \, \diffd W^k_v}^2 \giv
    \sF_{\sg_k}}
    \leq 16 c \int_0^\th \abs{a} (2\pi c v)^{-1/2} \, \dv
    \leq 16 \abs{a} (c\th)^{1/2} .
  \end{align*}
  Hence
  $\ee*{\xi_k \giv \sF_{\sg_k}} \leq (2\abs{a} + 4 \abs{a}^{1/2}
  (c\th)^{1/4})/c$, and the claim follows as in
  \ref{i-pieces-bd}.
\end{proof}

\section{Representation of the nonlinear law}\label{a-nonlinear-representation}\label{ap-thin}
We give the proof deferred from Section~\ref{ss-nlunique}, using
the notation of that section.
\begin{proof}[Proof of Lemma~\ref{l-thin}]
  We work on the canonical space
  $(\Om_1, \sF^P_\infty, \bF^P, P)$, where $\bF^P = (\sF^P_t)_{t \geq 0}$
  is the usual $P$-augmentation of the filtration $\bF^0$ generated
  by the coordinate process $(X, B, \ld)$, as in
  Section~\ref{ss-nldef}, and
  $\sF^P_\infty := \bigvee_t \sF^P_t$. The proof of
  Proposition~\ref{p-unique}, with one particle and $\bF^P$ in place
  of $\bF$, gives, on an extension of this space, a Poisson random
  measure $\Pi'$ with intensity $\ds \, \nu(\dz) \, \dy$ such that
  $(X_0, \ld_0)$, $B$ and $\Pi'$ are independent and
  $\mu(\ds,\dz) = \int \inds*{y \leq \ld_s} \, \Pi'(\ds,\dz,\dy)$. That
  proof uses only the fact that $\mu$ has unit atoms at distinct times,
  finitely many on bounded intervals, that its compensator is
  $\ld_s \, \ds \, \nu(\dz)$ with $\ld$ continuous and adapted, and that
  $B$ is a Brownian motion. By \ref{i-nl-graph} and \eqref{feed},
  $(X, \ld, \mu)$ then solves \eqref{frozen} with $v = u$ and $\Pi'$
  in place of $\Pi$, and
  $\int_0^T \norm*{u(s,\cdot)}_\infty \, \ds < \infty$ by
  Lemma~\ref{l-densbd}. By the uniqueness of solutions of
  \eqref{frozen},
  $(X, \ld, \mu) = \Phi_u(X_0, \ld_0, B, \Pi')$. The noise
  $(X_0, \ld_0, B, \Pi')$ has the same law as the noise
  $(X_0, \ld_0, B, \Pi)$ on $\Om$, and $\Phi_u$ is measurable.
  Hence $(\hX, B, \hat\ld)$ has the law of $(X, B, \ld)$ under $P$,
  which is $P$.
\end{proof}

\section{Measure-theoretic tools}\label{a-measure-tools}
We prove the two lemmas stated in Section~\ref{ss-iset}, in the
order in which they appear there.

\subsection{Extension from cylinder functionals}\label{ap-cylinder-extension}
\begin{proof}[Proof of Lemma~\ref{l-ext}]
  Put $\sT_s := \sT \cap [0,s]$, with $\sT_\infty := \sT$.
  Fix $m \in \N$ and $t_1 < \dots < t_m$ in $\sT_s$, and write
  $U := (Y_{t_1},\dots,Y_{t_m})$, which takes values in
  $H_m := (\R \times \R \times [0,\infty))^m$.
  The hypothesis gives $\ee*{\phi(U)\xi}=0$ for every
  $\phi \in \sC_m$. Since $\sC_m$ is dense in $C_c(H_m)$ for
  the supremum norm and
  \begin{align*}
    \abs*{\ee{(\phi(U)-\psi(U))\xi}}
    \leq \norm{\phi-\psi}_\infty\,\ee*{\abs{\xi}},
  \end{align*}
  the same identity holds for every $\phi \in C_c(H_m)$.
  The finite signed measure
  $\alpha(A):=\ee*{\xi\inds{U\in A}}$ on $H_m$ is determined
  by its integrals against $C_c(H_m)$, so $\alpha=0$.
  Thus $\ee*{\xi\inds{U\in A}}=0$ for every Borel set
  $A\subseteq H_m$, including $A=H_m$.

  The events of this form, over all finite sets of times in
  $\sT_s$, form a $\pi$-system generating
  $\sH_s:=\sigma(Y_t:t\in\sT_s)$. Since $\ee*{\xi}=0$,
  the events $C$ satisfying $\ee*{\xi\ind{C}}=0$ form a
  $\lambda$-system. The $\pi$--$\lambda$ theorem therefore gives
  $\ee*{\xi\giv\sH_s}=0$.

  Finally, $Y$ has right-continuous paths. For every finite
  $v<s$, choose $t_n\in\sT_s$ decreasing to $v$; then
  $Y_{t_n}\to Y_v$. If $s<\infty$, the endpoint $Y_s$ is
  already included because $s\in\sT$. This also covers $s=0$;
  when $s=\infty$, the same approximation applies to every
  finite $v$. Consequently,
  $\sH_s=\sigma(Y_v:v\leq s)$, where for $s=\infty$ only
  finite times $v$ are meant. This proves the claim.
\end{proof}

\subsection{Limits of random measures}\label{a-random-measures}
We give the proof of Lemma~\ref{l-pass}, using the notation of its statement.

\begin{proof}[Proof of Lemma~\ref{l-pass}]
  By monotone convergence, $\ee*{\inner*{M_N}{G}}$ $= \int G \, \diffd m_N$
  and $\int \ee^{\rho}{G} \, \mathfrak M(\diffd\rho) = \int G \, \diffd m$
  for every nonnegative
  measurable $G$. The set of discontinuity points of $F$ is
  $m$-null, so it is $\rho$-null for $\mathfrak M$-almost every
  $\rho$. For
  bounded such $F$, the map $\rho \mapsto \ee^{\rho}{F}$ on
  $\sP(S)$ is
  therefore continuous at $\mathfrak M$-almost every $\rho$, by
  the bounded, almost-everywhere continuous version of the
  Portmanteau theorem. The continuous mapping theorem then gives
  weak convergence of the laws of $\inner*{M_N}{F}$ to the image of
  $\mathfrak M$ under this map.

  \ref{i-pass-pos} For $K > 0$, this applies to $F \wedge K$, so
  $\int \ee^{\rho}{F \wedge K} \, \mathfrak M(\diffd\rho)$
  $= \lim_N \int (F \wedge K) \,
  \diffd m_N \leq \liminf_N \int F \, \diffd m_N$. Let
  $K \uparrow \infty$. A lower semicontinuous
  $F \colon S \to [0,\infty]$ is an increasing limit of bounded
  continuous functions, and the same argument applies to them.

  \ref{i-pass-zero} For $K > 0$, put $F_K := (F \wedge K) \vee (-K)$,
  so that $\abs{F - F_K} = (\abs{F} - K)^+$, and
  $\eps_K := \sup_N \int (\abs{F} - K)^+ \, \diffd m_N$, which tends to
  $0$ as $K \to \infty$. Then
  $\ee*{\abs*{\inner*{M_N}{F} - \inner*{M_N}{F_K}}}$ $\leq \eps_K$, and by
  \ref{i-pass-pos},
  $\int \ee^{\rho}{(\abs{F} - K)^+} \, \mathfrak M(\diffd\rho)
  \leq \eps_K$. In
  particular,
  $\int \ee^{\rho}{\abs{F}} \, \mathfrak M(\diffd\rho)
  \leq K + \eps_K$. Since
  the laws of $\inner*{M_N}{F_K}$ converge weakly to the image of
  $\mathfrak M$ under $\rho \mapsto \ee^{\rho}{F_K}$, with values in
  $[-K,K]$, we get, for every $K > 0$,
  \begin{align*}
    \int \abs*{\ee^{\rho}{F}} \, \mathfrak M(\diffd\rho)
    &\leq \textstyle\lim_N \ee{\abs{\inner{M_N}{F_K}}} + \eps_K
  \\ &\leq \textstyle \liminf_N \ee{\abs{\inner{M_N}{F}}} + 2 \eps_K
    = 2 \eps_K.  \qedhere
  \end{align*}
\end{proof}

\fin